\documentclass[11pt,a4paper,DIV=12]{scrartcl}

\usepackage[T1]{fontenc}
\usepackage{lmodern}
\usepackage{amsmath,amssymb,mathtools,amsthm}
\usepackage{dsfont}
\usepackage{stmaryrd}
\usepackage{enumitem}
\usepackage{microtype}
\usepackage[colorlinks,citecolor=blue,linkcolor=blue,urlcolor=blue]{hyperref}

\newtheoremstyle{upright}{\topsep}{\topsep}{\upshape}{}{\bfseries}{.}{ }{}
\theoremstyle{upright}
\newtheorem{theorem}{Theorem}[section]
\newtheorem{lemma}[theorem]{Lemma}
\newtheorem{proposition}[theorem]{Proposition}
\newtheorem{corollary}[theorem]{Corollary}
\theoremstyle{definition}
\newtheorem{definition}{Definition}[section]
\newtheorem{example}[definition]{Example}

\makeatletter
\renewenvironment{proof}[1][\proofname]{\par\pushQED{\qed}\normalfont
  \topsep6\p@\@plus6\p@\relax\trivlist
  \item[\hskip\labelsep\bfseries #1\@addpunct{.}]\ignorespaces}
  {\popQED\endtrivlist\@endpefalse}
\makeatother

\newcommand{\R}{\mathbf{R}}
\newcommand{\N}{\mathbf{N}}
\newcommand{\E}{\mathbf{E}}
\renewcommand{\P}{\mathbf{P}}
\newcommand{\Normal}{\mathcal{N}}

\newcommand{\ind}[1]{\mathds{1}_{\{ #1 \}}}
\newcommand{\id}{\operatorname{Id}}
\newcommand{\dom}[1]{\operatorname{Dom}{#1}}
\DeclareMathOperator{\var}{var}
\DeclareMathOperator{\cov}{cov}
\DeclareMathOperator{\oD}{d}
\newcommand{\dif}{\,\oD\!}
\newcommand{\sminus}{\setminus}
\newcommand{\bigCI}{\perp\mkern-10mu\perp}
\newcommand{\cyl}{\mathcal{S}}
\newcommand{\excond}[2]{\mathbf{E}\left[ \left.#1 \hspace{.2em}\right|\hspace{.2em} #2\right]}

\newcommand{\espxz}[1]{\mathbf{E}\left[#1 \mid X, Z \right]}

\newcommand{\deltaap}{\Delta^{\{a\}'}}

\title{Malliavin structure for conditionally independent random variables}
\author{Laurent Decreusefond\thanks{LTCI, T\'el\'ecom Paris, Institut
    Polytechnique de Paris, Palaiseau, France.\newline
    \texttt{laurent.decreusefond@mines-telecom.fr}}
  \and Christophe Vuong\thanks{T\'el\'ecom Paris, Palaiseau, France.\newline
    \texttt{christophe.vuong@telecom-paris.org}}}
\date{}

\begin{document}

\maketitle

\begin{abstract}
  We consider a countable family of random variables which are conditionally
  independent given a latent random variable $Z$. Given $Z$, the random variables are independent. Hence the Malliavin
  calculus for independent random variables applies for each value of $Z$. We show that these constructions fit
  together. We obtain the chaos decomposition, the domains of the Malliavin
  operators, the Mehler formula and the carré du champ operator, with
  conditional expectations given $Z$ in place of ordinary expectations. For normal
  approximation, conditioning on $Z$ and the infinite de Jong theorem of
  D\"obler give mixed normal approximations of degenerate U-statistics, in
  Wasserstein and Kolmogorov distances. Working under the unconditional
  probability, we obtain a de Jong type bound for U-statistics whose kernels
  are products. In both cases, the fluctuations of the conditional variance
  appear as an additional term, which cannot be removed. As an application, we
  consider subhypergraph counts in a random $3$-uniform hypergraph built on a
  latent random graph. For the count centred at its conditional mean, we
  obtain an explicit rate in Wasserstein distance. It reduces to the known
  rate for the Erd\H{o}s--R\'enyi hypergraph when the latent graph is
  complete. The count centred at its mean is a pattern count in a non-uniform
  hypergraph with independent edges, which gives a necessary and sufficient
  condition for its asymptotic normality.


\end{abstract}

\section{Introduction}
\label{sec:introduction}

Malliavin calculus is also known as the stochastic calculus of variations. At
the very core of it, it considers a gradient on a probability space. The link
between the gradient and the measure is made through an integration by parts
formula. On the Wiener space, the gradient generalizes the usual gradient on
$\R^{N}$ and the integration by parts yields an extension of the Itô integral.
On the space of point processes equipped with the law of a Poisson process, the
gradient becomes a difference operator and the integration by parts is an
avatar of the Mecke formula. For independent random variables, a Malliavin
structure was introduced in \cite{decreusefond:hal-01565240} and,
independently, in \cite{duerinckx2021size}; see also \cite{dung2018poisson}.
The gradient is then $D_{a}F = F - \E[F \mid X_{b},\, b \neq a]$. Recently,
\cite{dobler2026malliavin} completed this theory. The infinite Hoeffding
decomposition plays the role of the Wiener-Itô chaos decomposition, it gives
the domains of the Malliavin operators and the Mehler formula. Furthermore,
\cite{dobler2026malliavin} proves an infinite version of the quantitative de
Jong theorem.

By order of complexity, the next situation which can be analyzed is that of
conditionally independent random variables. This structure is very common. By
de Finetti's theorem, an infinite exchangeable sequence is conditionally
independent and identically distributed given its directing random measure.
Similarly, by the Aldous-Hoover theorem, the hyperedge indicators of an
exchangeable random hypergraph are conditionally independent given a latent
variable. This is the key to the analysis of random hypergraphs, see
\cite{austin_exchangeable_2008}.

Given the latent variable $Z = z$, the random variables are independent. Thus,
all the constructions of \cite{decreusefond:hal-01565240,dobler2026malliavin}
can be made for each value of $z$. The only point is to check that they fit
together, measurably in $z$. This is the content of
Section~\ref{sec:structure}. We obtain the chaos decomposition, the domains of
the gradient, of the divergence and of the number operator, the Mehler formula
and the carré du champ operator. At a formal level, expectations are replaced
by conditional expectations given $Z$. Hence, we give the statements and only
the arguments which are specific to the conditional setting. Note that the
kernel of the number operator is no longer the set of constants but the set of
$\sigma(Z)$-measurable random variables. Moreover, the gradient commutes with
the Ornstein-Uhlenbeck semigroup
as was observed in
\cite{dobler2026malliavin,duerinckx2021size} for independent random variables.

The prevailing application of Malliavin calculus is nowadays the evaluation of
convergence rates via Stein's method
\cite{NourdinNormalApproximationsMalliavin2012,Decreusefond2022}. In the
present setting, there are two routes. The first one is to condition on $Z$.
The bounds of \cite{dobler2026malliavin} then hold for each value of $Z$.
Their average bounds the distance between the law of a functional $F$ and a
variance mixture of Gaussian distributions. Furthermore, the distance between
this mixture and the standard Gaussian distribution is controlled by the
fluctuations of the conditional variance $\E[F^{2} \mid Z]$. This yields mixed
normal approximations of degenerate U-statistics in Wasserstein and Kolmogorov
distances (Theorem~\ref{thm-discrete-mall:mixed-dejong}). Note that the
resulting bounds are averages of ratios of conditional moments. The second
route works directly under the unconditional probability. We use Stein's method
and the carré du champ operator of the conditional structure, as in
\cite{decreusefond:hal-01565240,dobler2023normal}. The bounds are then
expressed in terms of unconditional moments. For U-statistics whose kernels
are products, we obtain a de Jong type bound with an explicit sum over
connected quadruples (Theorem~\ref{thm-discrete-mall:dejong-product}). In both
routes, the fluctuations of the conditional variance appear as an additional
term, which cannot be removed (Example~\ref{ex-discrete-mall:latent-term}).

Our motivation comes from random hypergraphs. We consider the layered model
$\mathbb{T}^{(3)}(n, q_{n}, p_{n})$. A triple carries a hyperedge only if its
three pairs are edges of a latent Erd\H{o}s-R\'enyi graph
$\mathbb{G}(n,q_{n})$, and then independently with probability $p_{n}$. Given
the latent graph, the hyperedges are independent. Let $N_{G}$ be the number of
copies of a $3$-uniform hypergraph $G$. For $N_{G}$ centred at its conditional
mean, we obtain an explicit rate in Wasserstein distance
(Theorem~\ref{thm-motif_estimation:dejong-quant-hypergraph}). When the latent
graph is complete, it reduces to the rate known for the Erd\H{o}s-R\'enyi
hypergraph. The count centred at its mean behaves differently. Its variance
splits into a hyperedge part and a latent part. Furthermore, $N_{G}$ is a
pattern count in a non-uniform random hypergraph with independent edges. Hence
the results of \cite{michalczuk2024normal} give a necessary and sufficient
condition for its asymptotic normality, with rates in Wasserstein and
Kolmogorov distances
(Theorem~\ref{thm-discrete-mall:hypergraph-condindep-threshold-hoeffding}).

Normal approximation for subgraph counts in random graphs has a long history.
It was first handled with the method of moments \cite{rucinski1988small}, which
yields the thresholds of asymptotic normality in terms of the number of
vertices but requires tedious computations. In \cite{barbour1989central},
Stein's method was used to derive convergence rates in the 1-Wasserstein
distance. The combination with Malliavin calculus has brought another feature
to the coupling constructions of Stein's method, through the chaos
representation property of Rademacher functionals \cite{privault2008stochastic},
with applications to subgraph counts \cite{privault2020normal} and to
percolation \cite{krokowski2017discrete}. For random hypergraphs whose
hyperedges are independent, asymptotic normality of subhypergraph counts goes
back to \cite{de1996central}. Necessary and sufficient conditions, with rates
in both the Wasserstein and the Kolmogorov distances, were recently obtained in
\cite{michalczuk2024normal}.

The paper is organized as follows. Section~\ref{sec:structure} builds the
Malliavin structure. Section~\ref{sec:normal-approximation} is devoted to
normal approximation, first by conditioning, then under the unconditional
probability. The application to random hypergraphs is in
Section~\ref{sec:appli}.

\section{Malliavin structure}
\label{sec:structure}

\subsection{Setting}
\label{sec:setting}

Let $A$ be an at most denumerable set. Let $E_{0}$ and $(E_{a},\, a\in A)$ be
Polish spaces and $E_{A} := \prod_{a\in A}E_{a}$, equipped with the product
Borel $\sigma$-field. On a probability space $(\Omega,\mathcal{T},\mathbb{P})$,
let $Z$ be an $E_{0}$-valued random variable and $X = (X_{a},\, a\in A)$ an
$E_{A}$-valued random variable. We assume that the $X_{a}$'s are conditionally
independent given $Z$:
\begin{equation*}
  X_a \operatornamewithlimits{\bigCI}_{Z} (X_b,\, b \in A \sminus \{a \}) \quad \text{for all } a\in A.
\end{equation*}
See chapter 5 of \cite{kallenberg1997foundations} for a review of conditional
independence. Since the spaces are Polish, $X_{a}$ admits a regular conditional
distribution given $Z$, denoted by $\mu_{a}(Z, \dif x_{a})$. By conditional
independence, the finite dimensional conditional distributions of $X$ given
$Z$ are products. Hence, for $\mathbb{P}_{Z}$-almost every $z$, the conditional
distribution of $X$ given $Z = z$ is the product measure
\begin{equation*}
  \P^{z} := \bigotimes_{a\in A}\mu_{a}(z, \cdot).
\end{equation*}
We denote by $\E^{z}$ the expectation with respect to $\P^{z}$.

For $p\ge 1$, $L^{p}(E_{A})$ is the space of random variables of the form
$F = f(Z, X)$, where $f\,:\, E_{0}\times E_{A}\to \R$ is measurable and
$\E[|F|^{p}] < \infty$. We identify $F$ and $f$, and we write $F(z, x)$. By
disintegration,
\begin{equation*}
  \excond{F}{Z} = \E^{Z}\bigl[F(Z, \cdot)\bigr] \quad\text{and}\quad \E[F] = \E\Bigl[\E^{Z}\bigl[F(Z,\cdot)\bigr]\Bigr].
\end{equation*}
In particular, if $F\in L^{p}(E_{A})$ then $F(z,\cdot)\in L^{p}(\P^{z})$ for
$\mathbb{P}_{Z}$-almost every $z$. A functional is cylindrical if it depends on
$X$ only through finitely many coordinates. By the martingale convergence
theorem, the set $\cyl$ of cylindrical functionals is dense in $L^{2}(E_{A})$.

For $x\in E_{A}$ and $K\subset A$, we write $x_{K} = (x_{a},\, a\in K)$ and
$x^{K} = (x_{a},\, a\in A\sminus K)$. We set
\begin{equation*}
  \mathcal{G}_{K} := \sigma(Z)\vee\sigma(X_{K}),\quad \mathcal{G}^{K} := \sigma(Z)\vee \sigma(X^{K}) \quad\text{and}\quad \mathcal{G}^{a} := \mathcal{G}^{\{a\}}.
\end{equation*}
By conditional independence,
\begin{equation}
  \label{eq-discrete-mall:cond-exp-kernel}
  \excond{F}{\mathcal{G}_{K}} = \int F(Z, X_{K}, y)\ \bigotimes_{b\notin K}\mu_{b}(Z, \dif y_{b}).
\end{equation}
The right-hand side is a jointly measurable function of $(Z, X_{K})$. For a
fixed $z$, it is the conditional expectation of $F(z,\cdot)$ given the
coordinates in $K$ under the product measure $\P^{z}$. Thus, for a fixed $z$,
the map $x\mapsto F(z,x)$ is a functional of independent random variables, in
the sense of \cite{decreusefond:hal-01565240,duerinckx2021size,dobler2026malliavin},
and \eqref{eq-discrete-mall:cond-exp-kernel} shows that the constructions of
these papers, made for each $z$, are jointly measurable in $(z,x)$.

We also need a conditionally independent copy of $X$. Up to an enlargement of
the probability space, let $X' = (X'_{a},\, a\in A)$ be such that, given $Z$,
$X$ and $X'$ are independent with the same conditional distribution
$\P^{Z}$. For $a\in A$, we set
\begin{equation*}
  F^{\{a\}'} := F(Z, X^{\{a\}}, X'_{a}) \quad\text{and}\quad \deltaap F := F - F^{\{a\}'}.
\end{equation*}
Given $\mathcal{G}^{a}$, the random variables $X_{a}$ and $X'_{a}$ are
independent with the same distribution. Hence $(F, F^{\{a\}'})$ and
$(F^{\{a\}'}, F)$ have the same distribution.

\subsection{Chaos decomposition}
\label{sec:chaos-decomp}

Let $\mathcal{P}_{f}(A)$ be the set of finite subsets of $A$. For
$F\in L^{2}(E_{A})$ and $M\in\mathcal{P}_{f}(A)$, set
\begin{equation*}
  F_{M} := \sum_{L\subset M}(-1)^{|M|-|L|}\,\excond{F}{\mathcal{G}_{L}}.
\end{equation*}
In particular, $F_{\emptyset} = \excond{F}{Z}$.

\begin{theorem}[Chaos decomposition]
  \label{thm-discrete-mall:chaos}
  Let $F\in L^{2}(E_{A})$.
  \begin{enumerate}
    \item For $M\in\mathcal{P}_{f}(A)$, $F_{M}$ is $\mathcal{G}_{M}$-measurable.
          For $K\subset A$, $\excond{F_{M}}{\mathcal{G}_{K}} = 0$ unless
          $M\subset K$.
    \item The $F_{M}$'s are pairwise orthogonal in $L^{2}(E_{A})$, and
          \begin{equation*}
            F = \sum_{M\in\mathcal{P}_{f}(A)}F_{M},
          \end{equation*}
          where the series converges unconditionally in $L^{2}(E_{A})$.
    \item If $F = \sum_{M}G_{M}$ in $L^{2}(E_{A})$, where the $G_{M}$'s satisfy
          item 1, then $G_{M} = F_{M}$ for every $M\in\mathcal{P}_{f}(A)$.
  \end{enumerate}
\end{theorem}

\begin{proof}
  By \eqref{eq-discrete-mall:cond-exp-kernel}, for
  $\mathbb{P}_{Z}$-almost every $z$, $F_{M}(z,\cdot)$ is the Hoeffding
  component of $F(z,\cdot)$ indexed by $M$ under $\P^{z}$, see
  \cite[Eq.~(2)]{dobler2026malliavin}. Moreover, $F_{M}$ is
  $\mathcal{G}_{M}$-measurable by construction. Hence, by
  \eqref{eq-discrete-mall:cond-exp-kernel},
  \begin{equation*}
    \excond{F_{M}}{\mathcal{G}_{K}} = \excond{F_{M}}{\mathcal{G}_{K\cap M}}.
  \end{equation*}
  Thus item 1 follows from the degeneracy of the Hoeffding components under
  each $\P^{z}$ \cite[Section~2.2]{dobler2026malliavin}.

  Let $M\neq N$, say $a\in M\sminus N$. Then $F_{N}$ is
  $\mathcal{G}^{a}$-measurable. Furthermore,
  $\mathcal{G}^{a} = \mathcal{G}_{A\sminus\{a\}}$, hence item 1 gives
  \begin{equation*}
    \E[F_{M}F_{N}] = \E\bigl[F_{N}\,\excond{F_{M}}{\mathcal{G}^{a}}\bigr] = 0.
  \end{equation*}
  By \cite[Proposition~2.1]{dobler2026malliavin}, for $\mathbb{P}_{Z}$-almost
  every $z$,
  \begin{equation*}
    \sum_{M\in\mathcal{P}_{f}(A)}\E^{z}[F_{M}^{2}] = \E^{z}[F^{2}] \quad\text{and}\quad \E^{z}[F F_{M}] = \E^{z}[F_{M}^{2}].
  \end{equation*}
  Let $(\mathcal{M}_{k})_{k}$ be an increasing sequence of finite subsets of
  $\mathcal{P}_{f}(A)$ whose union is $\mathcal{P}_{f}(A)$. Then
  \begin{equation*}
    \E\Bigl[\Bigl(F - \sum_{M\in\mathcal{M}_{k}}F_{M}\Bigr)^{2}\Bigr] = \E[F^{2}] - \sum_{M\in\mathcal{M}_{k}}\E[F_{M}^{2}],
  \end{equation*}
  which tends to $0$ by monotone convergence. Since the summands are orthogonal,
  the convergence is unconditional.

  For item 3, let $L\in\mathcal{P}_{f}(A)$. By item 1 and the continuity of
  the conditional expectation,
  \begin{equation*}
    \excond{F}{\mathcal{G}_{L}} = \sum_{N\subset L}G_{N}.
  \end{equation*}
  The Möbius inversion formula on the subsets of $M$ then gives $G_{M} = F_{M}$.
\end{proof}

For $p\ge 0$, set $\pi_{p}(F) := \sum_{|M| = p}F_{M}$ and let
$\mathfrak{C}_{p}$ be the range of $\pi_{p}$. By
Theorem~\ref{thm-discrete-mall:chaos}, the $\pi_{p}$'s are mutually orthogonal
projections and
\begin{equation*}
  F = \excond{F}{Z} + \sum_{p=1}^{\infty}\pi_{p}(F), \qquad
  L^{2}(E_{A}) = \bigoplus_{p=0}^{\infty}\mathfrak{C}_{p}.
\end{equation*}
Furthermore, $\mathfrak{C}_{0} = L^{2}(\Omega, \sigma(Z), \mathbb{P})$. For
$p\ge1$, $F$ belongs to $\mathfrak{C}_{p}$ if and only if
\begin{equation*}
  F = \sum_{M\subset A,\ |M| = p}W_{M},
\end{equation*}
where the series converges in $L^{2}(E_{A})$, each $W_{M}$ is
$\mathcal{G}_{M}$-measurable and $\excond{W_{M}}{\mathcal{G}^{a}} = 0$ for
every $a\in M$. Indeed, if $M\not\subset K$ and $a\in M\sminus K$, then
$\mathcal{G}_{K\cap M}\subset\mathcal{G}^{a}$, so that the last condition
implies item 1 of Theorem~\ref{thm-discrete-mall:chaos}. In statistical terms,
$F$ is a degenerate U-statistic of order $p$ whose kernels may depend on $Z$.
Under $\P^{z}$, $\mathfrak{C}_{p}$ is the $p$-th Hoeffding space of
\cite{dobler2026malliavin}.

\subsection{Malliavin operators}
\label{sec:mall-operators}

Let $L^{2}(A\times E_{A})$ be the Hilbert space of processes
$U = (U_{a},\, a\in A)$ such that
$\|U\|^{2}_{L^{2}(A\times E_{A})} := \sum_{a\in A}\E[U_{a}^{2}] < \infty$.

\begin{definition}[Discrete gradient]
  For $F\in L^{2}(E_{A})$ and $a\in A$, set
  \begin{equation*}
    D_{a}F := F - \excond{F}{\mathcal{G}^{a}}.
  \end{equation*}
  The domain $\dom D$ is the set of $F\in L^{2}(E_{A})$ such that
  $\sum_{a\in A}\E[(D_{a}F)^{2}] < \infty$, and then $DF := (D_{a}F,\, a\in A)\in L^{2}(A\times E_{A})$.
\end{definition}

\begin{theorem}
  \label{thm-discrete-mall:operators}
  Let $F\in L^{2}(E_{A})$.
  \begin{enumerate}
    \item For $a\in A$, $D_{a}F = \sum_{M\ni a}F_{M}$ and
          $D_{a}F = \espxz{\deltaap F}$. Each $D_{a}$ is an orthogonal
          projection, and $D_{a}D_{b} = D_{b}D_{a}$ for $a,b\in A$.
    \item One has
          \begin{equation*}
            \sum_{a\in A}\E[(D_{a}F)^{2}] = \sum_{M\in\mathcal{P}_{f}(A)}|M|\,\E[F_{M}^{2}] = \sum_{p=1}^{\infty}p\,\E[\pi_{p}(F)^{2}].
          \end{equation*}
          Moreover, $\dom D$ is the closure of $\cyl$ for the norm
          $\|F\|_{1,2}^{2} := \E[F^{2}] + \sum_{a}\E[(D_{a}F)^{2}]$.
    \item (Divergence) Let $U\in L^{2}(A\times E_{A})$. There exists $C$ such that
          $|\langle DG, U\rangle_{L^{2}(A\times E_{A})}| \le C\|G\|_{L^{2}(E_{A})}$
          for all $G\in\dom D$ if and only if the series
          $\sum_{a}D_{a}U_{a}$ converges in $L^{2}(E_{A})$. In this case,
          $\delta U := \sum_{a}D_{a}U_{a}$ satisfies
          \begin{equation*}
            \langle DG, U\rangle_{L^{2}(A\times E_{A})} = \E[G\,\delta U] \quad\text{for all } G\in\dom D.
          \end{equation*}
    \item (Number operator) Let $\dom{\mathsf{L}}$ be the set of
          $F\in L^{2}(E_{A})$ such that $\sum_{p}p^{2}\E[\pi_{p}(F)^{2}] < \infty$,
          and set $\mathsf{L}F := -\sum_{p}p\,\pi_{p}(F)$. Then
          $\mathsf{L}F = -\delta DF = -\sum_{a}D_{a}F$, the series converging in
          $L^{2}(E_{A})$. The operator $\mathsf{L}$ is self-adjoint, its spectrum
          is contained in $-\N_{0}$ and $\ker(\mathsf{L} + p\id) = \mathfrak{C}_{p}$.
          In particular, $\ker\mathsf{L} = L^{2}(\Omega,\sigma(Z),\mathbb{P})$.
    \item If $\excond{F}{Z} = 0$, then
          $\mathsf{L}^{-1}F := -\sum_{p\ge1}p^{-1}\pi_{p}(F)$ is the unique
          $G\in\dom{\mathsf{L}}$ such that $\excond{G}{Z} = 0$ and
          $\mathsf{L}G = F$.
  \end{enumerate}
\end{theorem}

\begin{proof}
  Let $M\in\mathcal{P}_{f}(A)$. Since $\mathcal{G}^{a} = \mathcal{G}_{A\sminus\{a\}}$,
  item 1 of Theorem~\ref{thm-discrete-mall:chaos} gives
  \begin{equation*}
    \excond{F_{M}}{\mathcal{G}^{a}} = \ind{a\notin M}F_{M}.
  \end{equation*}
  The conditional expectation is continuous on $L^{2}(E_{A})$, hence
  $D_{a}F = \sum_{M\ni a}F_{M}$. Furthermore, by
  \eqref{eq-discrete-mall:cond-exp-kernel},
  \begin{equation*}
    \espxz{F^{\{a\}'}} = \int F(Z, X^{\{a\}}, y)\,\mu_{a}(Z,\dif y) = \excond{F}{\mathcal{G}^{a}},
  \end{equation*}
  which gives the second expression of $D_{a}F$. The other assertions of item 1
  follow from the first expression. Item 2 is a consequence of item 1 and of
  the orthogonality of the $F_{M}$'s. The proofs of Propositions 2.5, 2.6, 2.8
  and 2.10 of \cite{dobler2026malliavin} only use the orthogonality of the
  Hoeffding components and the expression of $D_{a}F$ given in item 1. Thus
  they apply verbatim, and give the closure statement of item 2, item 3 and
  item 4. Item 5 follows from item 4.
\end{proof}

\subsection{Semigroup and carré du champ}
\label{sec:dirichlet}

Let $(\tau_{a},\, a\in A)$ be independent random variables with the exponential
distribution of mean $1$, independent of $(Z, X, X')$. For $t\ge 0$ and
$a, b\in A$, set
\begin{equation*}
  X_{b}(t) :=
  \begin{cases}
    X_{b}  & \text{if } \tau_{b} > t,   \\
    X'_{b} & \text{if } \tau_{b} \le t,
  \end{cases}
  \qquad
  X^{(a)}_{b}(t) :=
  \begin{cases}
    X_{b}    & \text{if } b = a,   \\
    X_{b}(t) & \text{if } b\neq a.
  \end{cases}
\end{equation*}
Given $Z$, $(X(t),\, t\ge0)$ is a Glauber dynamics: each coordinate $X_{b}$ is
refreshed at rate $1$, with a new value drawn from its conditional distribution
given $Z$. The process $(X^{(a)}(t),\, t\ge0)$ is the same dynamics in which
the coordinate $a$ is frozen.

\begin{theorem}
  \label{thm-discrete-mall:dirichlet}
  Let $F\in L^{2}(E_{A})$ and $t\ge0$.
  \begin{enumerate}
    \item (Mehler formula) Set $P_{t}F := \sum_{p\ge0}e^{-pt}\pi_{p}(F)$. Then
          \begin{equation*}
            P_{t}F = \excond{F(Z, X(t))}{Z, X}.
          \end{equation*}
          The family $(P_{t},\, t\ge0)$ is a strongly continuous contraction
          semigroup on $L^{2}(E_{A})$ whose generator is $\mathsf{L}$, and
          $P_{t}F$ tends to $\excond{F}{Z}$ in $L^{2}(E_{A})$ as $t\to\infty$.
          If $\excond{F}{Z} = 0$, then $-\mathsf{L}^{-1}F = \int_{0}^{\infty}P_{t}F\dif t$.
    \item (Commutation) For $a\in A$,
          \begin{equation*}
            D_{a}P_{t}F = P_{t}D_{a}F = e^{-t}P^{(a)}_{t}D_{a}F, \quad\text{where } P^{(a)}_{t}G := \excond{G(Z, X^{(a)}(t))}{Z, X}.
          \end{equation*}
    \item (Carré du champ) For $F, G\in\dom D$, set
          \begin{equation*}
            \Gamma(F,G) := \frac{1}{2}\sum_{a\in A}\espxz{\deltaap F\,\deltaap G}.
          \end{equation*}
          Then $\Gamma(F,G)\in L^{1}(E_{A})$ and
          $\excond{\Gamma(F,G)}{Z} = \sum_{a\in A}\excond{D_{a}F\,D_{a}G}{Z}$. If
          moreover $G\in\dom{\mathsf{L}}$, then
          \begin{equation*}
            \excond{F\,\mathsf{L}G}{Z} = -\excond{\Gamma(F,G)}{Z}.
          \end{equation*}
          If $F$, $G$ and $FG$ belong to $\dom{\mathsf{L}}$, then
          $2\Gamma(F,G) = \mathsf{L}(FG) - F\mathsf{L}G - G\mathsf{L}F$.
  \end{enumerate}
\end{theorem}

\begin{proof}
  For $\mathbb{P}_{Z}$-almost every $z$, given $Z = z$, the process $X(t)$ is
  the process of \cite[Section~2.1]{dobler2026malliavin} for the product
  measure $\P^{z}$, and $\pi_{p}(F)(z,\cdot)$ is the $p$-th component of the
  infinite Hoeffding decomposition of $F(z,\cdot)$. Thus the Mehler formula is
  \cite[Proposition~2.13]{dobler2026malliavin} under $\P^{z}$. The semigroup
  property, the contraction property, the limit as $t\to\infty$ and the
  identification of the generator follow from the spectral expression of
  $P_{t}$ and from Theorem~\ref{thm-discrete-mall:operators}, see
  \cite[Remark~2.12]{dobler2026malliavin}. The expression of $\mathsf{L}^{-1}$
  is \cite[Lemma~2.14]{dobler2026malliavin}, whose proof only uses the
  spectral expression.

  By item 1 of Theorem~\ref{thm-discrete-mall:operators},
  \begin{equation*}
    D_{a}P_{t}F = \sum_{M\ni a}e^{-t|M|}F_{M} = P_{t}D_{a}F.
  \end{equation*}
  Let $a\in M$ and write $F_{M} = f_{M}(Z, X_{M})$. Given $(Z, X)$, the
  coordinates $X^{(a)}_{b}(t)$, $b\in M\sminus\{a\}$, are independent, equal
  to $X_{b}$ with probability $e^{-t}$ and to $X'_{b}$ otherwise. If at least
  one of them is equal to $X'_{b}$, the conditional expectation of $f_{M}$
  vanishes by degeneracy. Hence, as in \cite[Eq.~(22)]{dobler2026malliavin},
  \begin{equation*}
    P^{(a)}_{t}F_{M} = e^{-t(|M|-1)}F_{M},
  \end{equation*}
  and item 2 follows.

  Item 3 is \cite[Eq.~(28), Propositions 2.17 and 2.18]{dobler2026malliavin}
  under $\P^{z}$.
\end{proof}

As in the case of independent random variables, functional inequalities follow
by conditioning.

\begin{corollary}
  \label{cor-discrete-mall:functional}
  Let $F\in L^{2}(E_{A})$.
  \begin{enumerate}
    \item (Conditional Efron-Stein inequality)
          \begin{equation*}
            \var(F\mid Z) \le \sum_{a\in A}\excond{(D_{a}F)^{2}}{Z} = \excond{\Gamma(F,F)}{Z}.
          \end{equation*}
          If $F\in\bigoplus_{k\le p}\mathfrak{C}_{k}$, the right-hand side is
          at most $p\,\var(F\mid Z)$.
    \item (Conditional McDiarmid inequality) Assume that there exist
          nonnegative numbers $(d_{a},\, a\in A)$ such that
          $\sum_{a}d_{a}^{2} < \infty$ and
          $|F(z, x) - F(z, x^{\{a\}}, y)| \le d_{a}$ for all $z$, $x$,
          $y\in E_{a}$ and $a\in A$. Then, for any $x>0$,
          \begin{equation*}
            \mathbb{P}\bigl(F - \excond{F}{Z} \ge x \mid Z\bigr) \le \exp\Bigl(-\frac{2x^{2}}{\sum_{a\in A}d_{a}^{2}}\Bigr).
          \end{equation*}
  \end{enumerate}
\end{corollary}

\begin{proof}
  Under $\P^{z}$, item 1 is the Efron-Stein inequality in the form of
  \cite[Eqs.~(29)-(31)]{dobler2026malliavin}. For item 2, let $(A_{m})_{m}$ be
  an increasing sequence of finite subsets of $A$ whose union is $A$, and
  \begin{equation*}
    F_{m} := \excond{F}{\mathcal{G}_{A_{m}}}.
  \end{equation*}
  By \eqref{eq-discrete-mall:cond-exp-kernel}, $F_{m}(z,\cdot)$ satisfies the
  same bounded differences condition. Thus, for every $z$, McDiarmid's
  inequality \cite{mcdiarmid1989method} applies to $F_{m}(z,\cdot)$ under the
  product measure $\P^{z}$. Since $\excond{F_{m}}{Z} = \excond{F}{Z}$, it gives,
  for $\varepsilon\in(0,x)$,
  \begin{equation*}
    \mathbb{P}\bigl(F_{m} - \excond{F}{Z} > x - \varepsilon \mid Z\bigr) \le \exp\Bigl(-\frac{2(x-\varepsilon)^{2}}{\sum_{a\in A}d_{a}^{2}}\Bigr).
  \end{equation*}
  Moreover, $F_{m}$ tends to $F$ almost surely, by the martingale convergence
  theorem. Hence the result follows by Fatou's lemma, letting $m\to\infty$ and
  then $\varepsilon\to0$.
\end{proof}

\section{Normal approximation}
\label{sec:normal-approximation}

For real random variables $U$ and $V$, the 1-Wasserstein and the Kolmogorov
distances are
\begin{equation*}
  d_{W}(U, V) := \sup_{h\in\textnormal{Lip}(1)}|\E[h(U)] - \E[h(V)]| \quad\text{and}\quad
  d_{K}(U,V) := \sup_{x\in\R}|\mathbb{P}(U\le x) - \mathbb{P}(V\le x)|,
\end{equation*}
where $\textnormal{Lip}(1)$ is the set of $1$-Lipschitz functions. We use the
same notation for probability distributions. In what follows, $\mathcal{N}$
denotes a standard Gaussian random variable, independent of $(Z, X, X')$, and
$\Phi$ its distribution function.

\subsection{Normal approximation by conditioning}
\label{sec:conditioning}

Let $F\in L^{2}(E_{A})$. We denote by $\mathcal{L}(F\mid Z)$ the conditional
distribution of $F$ given $Z$ and set $\sigma(Z) := \excond{F^{2}}{Z}^{1/2}$.
Let
\begin{equation*}
  \delta_{W}(Z) := d_{W}\bigl(\mathcal{L}(F\mid Z), \Normal(0, \sigma^{2}(Z))\bigr)
  \quad\text{and}\quad
  \delta_{K}(Z) := d_{K}\bigl(\mathcal{L}(F\mid Z), \Normal(0, \sigma^{2}(Z))\bigr),
\end{equation*}
with the convention that $\Normal(0,0)$ is the Dirac mass at $0$.

\begin{lemma}
  \label{lem-discrete-mall:mixture}
  Let $F\in L^{2}(E_{A})$.
  \begin{enumerate}
    \item For every measurable $h\,:\,\R\times E_{0}\to\R$ such that $h(\cdot, z)$
          is $1$-Lipschitz for every $z$,
          \begin{equation*}
            \bigl|\E[h(F, Z)] - \E[h(\sigma(Z)\mathcal{N}, Z)]\bigr| \le \E[\delta_{W}(Z)].
          \end{equation*}
          For every $x\in\R$ and every Borel set $B\subset E_{0}$,
          \begin{equation*}
            \bigl|\mathbb{P}(F\le x,\ Z\in B) - \mathbb{P}(\sigma(Z)\mathcal{N}\le x,\ Z\in B)\bigr| \le \E[\delta_{K}(Z)].
          \end{equation*}
    \item Furthermore, we have
          \begin{align*}
            d_{W}(F, \Normal(0,1)) & \le \E[\delta_{W}(Z)] + \sqrt{\frac{2}{\pi}}\ \E\bigl|\sigma(Z) - 1\bigr|, \\
            d_{K}(F, \Normal(0,1)) & \le \E[\delta_{K}(Z)] + \E\bigl[\ell(\sigma^{2}(Z))\bigr],
          \end{align*}
          where $\ell(v) := \min\bigl(1/2,\ |\log v|/(2\sqrt{2\pi e})\bigr)$ for
          $v>0$ and $\ell(0) := 1/2$.
  \end{enumerate}
\end{lemma}

\begin{proof}
  Since $\mathcal{N}$ is independent of $Z$, $\Normal(0,\sigma^{2}(Z))$ is the
  conditional distribution of $\sigma(Z)\mathcal{N}$ given $Z$. Hence
  \begin{equation*}
    \bigl|\excond{h(F,Z)}{Z} - \excond{h(\sigma(Z)\mathcal{N}, Z)}{Z}\bigr| \le \delta_{W}(Z),
  \end{equation*}
  and item 1 follows by integration. The same argument applies to the
  indicator functions of $(-\infty, x]\times B$. For item 2, the triangle
  inequality and item 1 reduce the problem to the distance between
  $\sigma(Z)\mathcal{N}$ and $\mathcal{N}$. On the one hand,
  \begin{equation*}
    \E\bigl|\sigma(Z)\mathcal{N} - \mathcal{N}\bigr| = \sqrt{\frac{2}{\pi}}\ \E\bigl|\sigma(Z) - 1\bigr|.
  \end{equation*}
  On the other hand, with the convention $\Phi(x/0) := \ind{x\ge0}$,
  \begin{equation*}
    \mathbb{P}(\sigma(Z)\mathcal{N}\le x) - \Phi(x) = \E\bigl[\Phi(x/\sigma(Z)) - \Phi(x)\bigr].
  \end{equation*}
  For $v>0$, set $u = xv^{-1/2}$. Then
  \begin{equation*}
    \frac{\dif}{\dif v}\,\Phi\bigl(xv^{-1/2}\bigr) = -\frac{u\,\Phi'(u)}{2v} \quad\text{and}\quad |u\Phi'(u)|\le\frac{1}{\sqrt{2\pi e}}.
  \end{equation*}
  Thus, for $\sigma > 0$,
  \begin{equation*}
    \bigl|\Phi(x/\sigma) - \Phi(x)\bigr| \le \frac{|\log\sigma^{2}|}{2\sqrt{2\pi e}}.
  \end{equation*}
  Moreover, $\Phi(x/\sigma)$ and $\Phi(x)$ both belong to $[0,1/2]$ if $x<0$
  and to $[1/2, 1]$ if $x\ge0$, so that their difference is at most $1/2$.
\end{proof}

Item 1 quantifies the convergence in distribution of $(F, Z)$ to
$(\sigma(Z)\mathcal{N}, Z)$. Thus, for a sequence $(F_{n})_{n}$ of functionals
of the same pair $(Z, X)$ such that $\E[\delta_{W,n}(Z)]$ tends to $0$ and
$\sigma_{n}(Z)$ converges in probability to some $\sigma(Z)$, the pair
$(F_{n}, Z)$ converges in distribution to $(\sigma(Z)\mathcal{N}, Z)$. In
other words, $F_{n}$ converges $\sigma(Z)$-stably to the mixture
$\sigma(Z)\mathcal{N}$. For degenerate
U-statistics, $\delta_{W}$ and $\delta_{K}$ are given by the infinite de Jong
theorem of \cite{dobler2026malliavin}. For $F\in L^{2}(E_{A})$, the
conditional maximal influence is
\begin{equation*}
  \rho_{Z}(F)^{2} := \sup_{a\in A}\excond{(D_{a}F)^{2}}{Z},
\end{equation*}
see \cite{dobler2017quantitative,mossel2010noise} for the role of this quantity.

\begin{theorem}
  \label{thm-discrete-mall:mixed-dejong}
  Let $p\ge1$ and $F\in\mathfrak{C}_{p}\cap L^{4}(E_{A})$. Set
  $\kappa(Z) := \excond{F^{4}}{Z} - 3\sigma^{4}(Z)$. With the convention that
  the ratios below vanish on $\{\sigma(Z) = 0\}$,
  \begin{align*}
    \E[\delta_{W}(Z)] & \le C_{p}\ \E\Bigl[\frac{\sqrt{|\kappa(Z)|}}{\sigma(Z)} + \rho_{Z}(F)\Bigr],                                                                                 \\
    \E[\delta_{K}(Z)] & \le \E\Bigl[\min\Bigl(1,\ 11.9\,\frac{\sqrt{|\kappa(Z)|}}{\sigma^{2}(Z)} + \bigl(3.5 + 10.8\sqrt{\kappa_{p}}\bigr)\frac{\rho_{Z}(F)}{\sigma(Z)}\Bigr)\Bigr],
  \end{align*}
  where $\kappa_{p}$ is the combinatorial constant of
  \cite{dobler2017quantitative} and $C_{p}$ only depends on $p$.
\end{theorem}

\begin{proof}
  Since $\pi_{q}(F) = 0$ for $q\neq p$, for $\mathbb{P}_{Z}$-almost every $z$,
  $F(z,\cdot)$ is a completely degenerate U-statistic of order $p$ under
  $\P^{z}$ and belongs to $L^{4}(\P^{z})$. If $\sigma(z) = 0$, then
  $F(z,\cdot) = 0$ $\P^{z}$-almost surely, and $\delta_{W}(z) = \delta_{K}(z) = 0$.
  Otherwise, set $G := F(z,\cdot)/\sigma(z)$. Then
  \begin{equation*}
    \E^{z}[G^{2}] = 1 \quad\text{and}\quad \E^{z}[G^{4}] - 3 = \frac{\kappa(z)}{\sigma^{4}(z)}.
  \end{equation*}
  Furthermore, by \eqref{eq-discrete-mall:cond-exp-kernel}, the influence of
  the coordinate $a$ on $G$ under $\P^{z}$ is
  \begin{equation*}
    \frac{\E\bigl[(D_{a}F)^{2}\mid Z = z\bigr]}{\sigma^{2}(z)},
  \end{equation*}
  so that the maximal influence of $G$ is $\rho_{z}(F)/\sigma(z)$. Thus
  \cite[Theorem~3.1]{dobler2026malliavin} bounds $d_{W}(G, \mathcal{N})$ and
  $d_{K}(G,\mathcal{N})$. Since
  \begin{equation*}
    \delta_{W}(z) = \sigma(z)\,d_{W}(G,\mathcal{N}) \quad\text{and}\quad \delta_{K}(z) = d_{K}(G,\mathcal{N})\le 1,
  \end{equation*}
  the result follows by integration.
\end{proof}

When $\E[F^{2}] = 1$, one has $\E[\sigma^{2}(Z)] = 1$, so that
\begin{equation*}
  \E\bigl|\sigma(Z) - 1\bigr| \le \E\bigl|\sigma^{2}(Z) - 1\bigr| \le \var\bigl[\sigma^{2}(Z)\bigr]^{1/2}.
\end{equation*}
Hence Lemma~\ref{lem-discrete-mall:mixture} and
Theorem~\ref{thm-discrete-mall:mixed-dejong} give
\begin{equation*}
  d_{W}(F, \Normal(0,1)) \le C_{p}\ \E\Bigl[\frac{\sqrt{|\kappa(Z)|}}{\sigma(Z)} + \rho_{Z}(F)\Bigr] + \sqrt{\frac{2}{\pi}}\ \var\bigl[\excond{F^{2}}{Z}\bigr]^{1/2}.
\end{equation*}
The last term cannot be removed.

\begin{example}
  \label{ex-discrete-mall:latent-term}
  Let $(\varepsilon_{i})_{i\ge1}$ be independent Rademacher random variables,
  independent of $Z$. Let $c(Z)$ take the values $\sqrt{1/2}$ and
  $\sqrt{3/2}$ with probability $1/2$, and
  $F = n^{-1/2}c(Z)\sum_{i=1}^{n}\varepsilon_{i}$. Then
  $F\in\mathfrak{C}_{1}$, $\E[F^{2}] = 1$, $\kappa(Z) = -2c(Z)^{4}/n$ and
  $\rho_{Z}(F) = c(Z)/\sqrt{n}$. Thus the first term of the previous bound is
  of order $n^{-1/2}$. However, $F$ converges in distribution to the mixture
  $c(Z)\mathcal{N}$, which is not Gaussian. Here
  $\var[\excond{F^{2}}{Z}] = 1/4$.
\end{example}

\subsection{Bounds under the unconditional probability}
\label{sec:unconditional}

The bounds of Theorem~\ref{thm-discrete-mall:mixed-dejong} are averages of
ratios of conditional moments. For a functional normalised by its unconditional
variance, it is often simpler to estimate unconditional moments. Thus we now
work directly under $\mathbb{P}$. We recall a standard implementation of
Stein's method for the Wasserstein distance to the Gaussian distribution, see
\cite{Decreusefond2022}. The bound on $\varphi''$ is
\cite[Lemma~2.4]{chen2010normal}.

\begin{lemma}
  \label{lem-discrete-mall:stein-normal}
  Let $\mathcal{H}_{*} := \{\varphi\in C^{2}(\R,\R)\,:\, \|\varphi'\|_{\infty}\le\sqrt{2/\pi},\ \|\varphi''\|_{\infty}\le 2\}$.
  For $F\in L^{1}(E_{A})$,
  \begin{equation*}
    d_{W}(F, \Normal(0,1)) \le \sup_{\varphi\in\mathcal{H}_{*}}\bigl|\E[\varphi'(F) - F\varphi(F)]\bigr|.
  \end{equation*}
\end{lemma}

The operator $\mathsf{L}$ is not diffusive: in general
$\Gamma(\psi(F), G)\neq\psi'(F)\Gamma(F,G)$. The following pseudo chain rule
measures the defect. Its statement and its proof follow the argument of
\cite[p.~177]{dobler2023normal}, to which we refer for the corresponding
computation in the independent setting.

\begin{lemma}[Pseudo chain rule]
  \label{lem-discrete-mall:approx-diffusion}
  Let $\psi\in C^{2}(\R,\R)$ with $\psi'$ and $\psi''$ bounded, and
  $F, G\in\dom D$. Then $\psi(F)\in\dom D$ and
  \begin{equation*}
    \Gamma(\psi(F), G) = \psi'(F)\,\Gamma(F,G) + R_{\psi}(F,G),
  \end{equation*}
  where
  \begin{equation*}
    |R_{\psi}(F,G)| \le \frac{\|\psi''\|_{\infty}}{4}\sum_{a\in A}\espxz{|\deltaap G|\,(\deltaap F)^{2}}.
  \end{equation*}
\end{lemma}

\begin{proof}
  By item 1 of Theorem~\ref{thm-discrete-mall:operators},
  \begin{equation*}
    \E\bigl[F F^{\{a\}'}\bigr] = \E\bigl[\excond{F}{\mathcal{G}^{a}}^{2}\bigr].
  \end{equation*}
  Since $F^{\{a\}'}$ has the same distribution as $F$, it follows that
  \begin{equation}
    \label{eq-discrete-mall:delta-D}
    \E\bigl[(\deltaap F)^{2}\bigr] = 2\E\bigl[(D_{a}F)^{2}\bigr].
  \end{equation}
  Moreover, $|\deltaap\psi(F)|\le\|\psi'\|_{\infty}|\deltaap F|$, hence
  $\psi(F)\in\dom D$. Furthermore,
  \begin{equation*}
    \deltaap\psi(F) = \psi'(F)\,\deltaap F - r_{a},
  \end{equation*}
  where
  \begin{equation*}
    r_{a} := \psi(F^{\{a\}'}) - \psi(F) - \psi'(F)\bigl(F^{\{a\}'} - F\bigr)
    \quad\text{satisfies}\quad
    |r_{a}|\le\frac{\|\psi''\|_{\infty}}{2}\,(\deltaap F)^{2}
  \end{equation*}
  by Taylor's formula. Since $\psi'(F)$ is $\sigma(Z, X)$-measurable, the
  result follows from the definition of $\Gamma$.
\end{proof}

\begin{theorem}
  \label{thm-discrete-mall:bounds-wass-general}
  Let $F\in\dom D$ such that $\excond{F}{Z} = 0$ and $\E[F^{2}] = 1$. Then
  \begin{equation*}
    d_{W}(F, \Normal(0,1)) \le \sqrt{\frac{2}{\pi}}\ \E\bigl|1 - \Gamma(F, -\mathsf{L}^{-1}F)\bigr| + \frac{1}{2}\sum_{a\in A}\E\bigl[|\deltaap\mathsf{L}^{-1}F|\,(\deltaap F)^{2}\bigr].
  \end{equation*}
  Moreover, if $F\in L^{4}(E_{A})$,
  \begin{multline}
    \label{eq-discrete-mall:bound-thm-1-Wasserstein}
    d_{W}(F, \Normal(0,1)) \le \sqrt{\frac{2}{\pi}}\ \var\bigl[\Gamma(F, -\mathsf{L}^{-1}F)\bigr]^{1/2}\\
    + \frac{\sqrt{2}}{2}\ \bigl(-\E[F\mathsf{L}^{-1}F]\bigr)^{1/2}\Bigl(\sum_{a\in A}\E\bigl[(\deltaap F)^{4}\bigr]\Bigr)^{1/2}.
  \end{multline}
\end{theorem}

\begin{proof}
  Let $\varphi\in\mathcal{H}_{*}$ and $G := -\mathsf{L}^{-1}F$, which belongs to
  $\dom{\mathsf{L}}$. Since $\excond{F}{Z} = 0$,
  \begin{equation*}
    F = -\mathsf{L}G.
  \end{equation*}
  Hence, by item 3 of Theorem~\ref{thm-discrete-mall:dirichlet} and
  Lemma~\ref{lem-discrete-mall:approx-diffusion},
  \begin{equation*}
    \E[F\varphi(F)] = -\E[\varphi(F)\,\mathsf{L}G] = \E[\Gamma(\varphi(F), G)] = \E[\varphi'(F)\Gamma(F,G)] + \E[R_{\varphi}(F,G)].
  \end{equation*}
  Thus
  \begin{equation*}
    \bigl|\E[\varphi'(F) - F\varphi(F)]\bigr| \le \|\varphi'\|_{\infty}\,\E|1 - \Gamma(F,G)| + \frac{\|\varphi''\|_{\infty}}{4}\sum_{a\in A}\E\bigl[|\deltaap G|\,(\deltaap F)^{2}\bigr],
  \end{equation*}
  and Lemma~\ref{lem-discrete-mall:stein-normal} gives the first bound. For the
  second one, note that
  \begin{equation*}
    \E[\Gamma(F,G)] = -\E[F\mathsf{L}G] = \E[F^{2}] = 1,
  \end{equation*}
  so that the Cauchy-Schwarz inequality bounds the first term. Furthermore, by
  \eqref{eq-discrete-mall:delta-D},
  \begin{equation*}
    \sum_{a\in A}\E\bigl[(\deltaap G)^{2}\bigr] = 2\sum_{a\in A}\E\bigl[(D_{a}G)^{2}\bigr] = 2\E[\Gamma(G, G)] = -2\E[F\mathsf{L}^{-1}F],
  \end{equation*}
  and the Cauchy-Schwarz inequality bounds the second term.
\end{proof}

\begin{corollary}
  \label{cor-discrete-mall:bound-multi-chaos}
  Let $F = \sum_{p=1}^{m}F_{p}\in L^{4}(E_{A})$ with $F_{p}\in\mathfrak{C}_{p}$
  and $\E[F^{2}] = 1$. Then
  \begin{multline*}
    d_{W}(F, \Normal(0,1)) \le \sqrt{\frac{2}{\pi}}\sum_{p,q=1}^{m}\frac{1}{q}\,\var\bigl[\Gamma(F_{p},F_{q})\bigr]^{1/2}\\
    + \frac{\sqrt{2}}{2}\Bigl(\sum_{p=1}^{m}\frac{\E[F_{p}^{2}]}{p}\Bigr)^{1/2}\Bigl\{\sum_{p=1}^{m}\Bigl(\sum_{a\in A}\E\bigl[(\deltaap F_{p})^{4}\bigr]\Bigr)^{1/4}\Bigr\}^{2}.
  \end{multline*}
\end{corollary}

\begin{proof}
  Since $-\mathsf{L}^{-1}F = \sum_{q}F_{q}/q$, one has
  \begin{equation*}
    \Gamma(F, -\mathsf{L}^{-1}F) = \sum_{p,q=1}^{m}\frac{1}{q}\,\Gamma(F_{p},F_{q}).
  \end{equation*}
  Thus the triangle inequality in $L^{2}(E_{A})$ bounds the first term of
  \eqref{eq-discrete-mall:bound-thm-1-Wasserstein}. For the second one, by
  orthogonality,
  \begin{equation*}
    -\E[F\mathsf{L}^{-1}F] = \sum_{p=1}^{m}\frac{\E[F_{p}^{2}]}{p}.
  \end{equation*}
  Moreover, Minkowski's inequality gives
  \begin{equation*}
    \Bigl(\sum_{a\in A}\E\bigl[(\deltaap F)^{4}\bigr]\Bigr)^{1/4} \le \sum_{p=1}^{m}\Bigl(\sum_{a\in A}\E\bigl[(\deltaap F_{p})^{4}\bigr]\Bigr)^{1/4}.
  \end{equation*}
\end{proof}

We now specialise to U-statistics whose kernels are products. In the rest of
this section, $F = \sum_{p=1}^{m}F_{p}$ with $F_{p} = \sum_{|I| = p}W_{I}$, a
finite sum, and
\begin{equation}
  \label{eq-discrete-mall:product-form}
  W_{I} = c_{I}\prod_{\alpha\in I}Y_{\alpha},
\end{equation}
where the $c_{I}$'s are $\sigma(Z)$-measurable with $\E[c_{I}^{4}] < \infty$.
Each $Y_{\alpha}$ is $\sigma(Z, X_{\alpha})$-measurable, conditionally centred
given $Z$ and bounded by $1$. In particular, the $Y_{\alpha}$'s are
conditionally independent given $Z$, $F_{p}\in\mathfrak{C}_{p}$ and, for every
integer $k\ge2$,
$\excond{|Y_{\alpha}|^{k}}{Z}\le\excond{Y_{\alpha}^{2}}{Z}\le1$.

\begin{definition}[Connected and admissible quadruples]
  \label{def-discrete-mall:connected-admissible}
  A quadruple $(I_{1}, \ldots, I_{4})$ of nonempty subsets of $A$ is connected
  if the graph with vertex set $\{I_{1}, \ldots, I_{4}\}$, in which $I_{i}$ and
  $I_{j}$ are adjacent when $I_{i}\cap I_{j}\neq\emptyset$, is connected. It is
  admissible if every element of $I_{1}\cup\dots\cup I_{4}$ belongs to at least
  two of the $I_{j}$'s.
\end{definition}

Let $\mathcal{Q}$ be the set of quadruples of index sets of kernels of $F$
which are admissible and connected. Set
\begin{equation}
  \label{eq-discrete-mall:Mplus}
  \mathfrak{M}^{+}(F) := \sum_{(I,J,K,L)\in\mathcal{Q}}\E\bigl[B(I,J,K,L)\bigr],
\end{equation}
where
\begin{equation*}
  B(I,J,K,L) := |c_{I}c_{J}c_{K}c_{L}|\prod_{\alpha\in I\cup J\cup K\cup L}\excond{Y_{\alpha}^{2}}{Z}.
\end{equation*}
Since the $Y_{\alpha}$'s are conditionally independent and conditionally
centred, a conditional expectation of the form
$\excond{\prod_{\alpha}Y_{\alpha}^{m_{\alpha}}}{Z}$ vanishes as soon as one of
the multiplicities $m_{\alpha}$ is equal to $1$. Otherwise, it is bounded by
$\prod_{\alpha}\excond{Y_{\alpha}^{2}}{Z}$, the product running over the
$\alpha$'s with $m_{\alpha}\ge 2$. This remark is used repeatedly below.

\begin{proposition}
  \label{prop-discrete-mall:var-gamma-product}
  For $1\le p,q\le m$,
  \begin{equation*}
    \var\bigl[\Gamma(F_{p},F_{q})\bigr] \le 2m^{2}\,\mathfrak{M}^{+}(F) + \ind{p=q}\,p^{2}\,\var\bigl[\excond{F_{p}^{2}}{Z}\bigr].
  \end{equation*}
\end{proposition}

\begin{proof}
  By item 3 of Theorem~\ref{thm-discrete-mall:dirichlet} and item 1 of
  Theorem~\ref{thm-discrete-mall:operators},
  \begin{equation*}
    \excond{\Gamma(F_{p},F_{q})}{Z} = \sum_{a\in A}\ \sum_{I\ni a,\ J\ni a}\excond{W_{I}W_{J}}{Z} = \ind{p=q}\,p\,\excond{F_{p}^{2}}{Z},
  \end{equation*}
  since $\excond{W_{I}W_{J}}{Z} = 0$ for $I\neq J$. Therefore
  \begin{equation*}
    \var\bigl[\Gamma(F_{p},F_{q})\bigr] = \E\bigl[\var(\Gamma(F_{p},F_{q})\mid Z)\bigr] + \ind{p=q}\,p^{2}\,\var\bigl[\excond{F_{p}^{2}}{Z}\bigr].
  \end{equation*}
  Let $a\in A$. Then
  \begin{equation*}
    \deltaap F_{p} = \sum_{I\ni a}\bigl(W_{I} - W_{I}^{\{a\}'}\bigr).
  \end{equation*}
  For $I, J\ni a$, by \eqref{eq-discrete-mall:cond-exp-kernel} and
  degeneracy,
  \begin{equation*}
    \espxz{W_{J}^{\{a\}'}} = \excond{W_{J}}{\mathcal{G}^{a}} = 0
    \quad\text{and}\quad
    \espxz{W_{I}^{\{a\}'}W_{J}^{\{a\}'}} = \excond{W_{I}W_{J}}{\mathcal{G}^{a}}.
  \end{equation*}
  Hence
  \begin{equation*}
    \espxz{(W_{I} - W_{I}^{\{a\}'})(W_{J} - W_{J}^{\{a\}'})} = W_{I}W_{J} + \excond{W_{I}W_{J}}{\mathcal{G}^{a}}.
  \end{equation*}
  Thus, with $E_{a} := \excond{\cdot}{\mathcal{G}^{a}}$,
  \begin{equation*}
    2\Gamma(F_{p}, F_{q}) = \sum_{(I,J,a)}T_{IJa}, \qquad T_{IJa} := W_{I}W_{J} + E_{a}(W_{I}W_{J}),
  \end{equation*}
  the sum running over $|I| = p$, $|J| = q$ and $a\in I\cap J$. By
  \eqref{eq-discrete-mall:product-form} and conditional independence,
  \begin{equation*}
    E_{a}(W_{I}W_{J}) = c_{I}c_{J}\,\excond{Y_{a}^{2}}{Z}\prod_{\alpha\in I\sminus\{a\}}Y_{\alpha}\prod_{\alpha\in J\sminus\{a\}}Y_{\alpha}.
  \end{equation*}
  Hence $T_{IJa}$ is measurable with respect to the $\sigma$-field
  \begin{equation*}
    \sigma(Z)\vee\sigma(Y_{\alpha},\, \alpha\in I\cup J).
  \end{equation*}
  It follows that $\cov(T_{IJa}, T_{KLb}\mid Z) = 0$ when
  $(I\cup J)\cap(K\cup L) = \emptyset$, and
  \begin{equation*}
    4\var\bigl(\Gamma(F_{p},F_{q})\mid Z\bigr) \le \sum_{(I\cup J)\cap(K\cup L)\neq\emptyset}\Bigl(\bigl|\excond{T_{IJa}T_{KLb}}{Z}\bigr| + \bigl|\excond{T_{IJa}}{Z}\bigr|\,\bigl|\excond{T_{KLb}}{Z}\bigr|\Bigr).
  \end{equation*}
  When $(I\cup J)\cap(K\cup L)\neq\emptyset$, the quadruple $(I,J,K,L)$ is
  connected, since $a\in I\cap J$ and $b\in K\cap L$.

  Consider $\excond{E_{a}(W_{I}W_{J})\,E_{b}(W_{K}W_{L})}{Z}$. The three other
  products in $\excond{T_{IJa}T_{KLb}}{Z}$ are treated in the same way. It is
  equal to
  \begin{equation*}
    c_{I}c_{J}c_{K}c_{L}\,\excond{Y_{a}^{2}}{Z}\excond{Y_{b}^{2}}{Z}\prod_{\alpha}\excond{Y_{\alpha}^{m'_{\alpha}}}{Z},
  \end{equation*}
  where $m'_{\alpha}$ is the multiplicity of $\alpha$ in
  $(I\sminus\{a\}, J\sminus\{a\}, K\sminus\{b\}, L\sminus\{b\})$. For
  $\alpha\notin\{a,b\}$, $m'_{\alpha}$ is the multiplicity of $\alpha$ in
  $(I,J,K,L)$. Moreover, $a$ and $b$ belong to at least two of the four sets.
  Hence the term vanishes unless $(I,J,K,L)$ is admissible, and it is then
  bounded by $B(I,J,K,L)$. Next,
  \begin{equation*}
    \excond{T_{IJa}}{Z} = 2\excond{W_{I}W_{J}}{Z}
  \end{equation*}
  vanishes unless $I = J$, and similarly for $(K, L)$. For $I = J$, $K = L$
  and $I\cap K\neq\emptyset$, the quadruple $(I,I,K,K)$ is admissible and
  connected, and
  \begin{equation*}
    \bigl|\excond{T_{IIa}}{Z}\,\excond{T_{KKb}}{Z}\bigr|\le 4B(I,I,K,K).
  \end{equation*}
  Thus each pair $((I,J,a),(K,L,b))$ contributes at most $8B(I,J,K,L)$ with
  $(I,J,K,L)\in\mathcal{Q}$. Each quadruple arises from at most $m^{2}$ pairs
  $(a,b)$, so that
  \begin{equation*}
    4\var\bigl(\Gamma(F_{p},F_{q})\mid Z\bigr)\le 8m^{2}\sum_{(I,J,K,L)\in\mathcal{Q}}B(I,J,K,L).
  \end{equation*}
  Taking expectations gives the result.
\end{proof}

\begin{lemma}
  \label{lem-discrete-mall:fourth-product}
  For $1\le p\le m$,
  $\sum_{a\in A}\E\bigl[(\deltaap F_{p})^{4}\bigr]\le 16\,p\,\mathfrak{M}^{+}(F)$.
\end{lemma}

\begin{proof}
  Let $a\in A$ and
  \begin{equation*}
    S_{a} := \sum_{I\ni a}W_{I},
  \end{equation*}
  the sum running over the kernels of $F_{p}$. Then
  $\deltaap F_{p} = S_{a} - S_{a}^{\{a\}'}$. Since
  $(x+y)^{4}\le 8(x^{4} + y^{4})$ and $S_{a}^{\{a\}'}$ has the same distribution
  as $S_{a}$,
  \begin{equation*}
    \E\bigl[(\deltaap F_{p})^{4}\bigr]\le 16\,\E\bigl[S_{a}^{4}\bigr].
  \end{equation*}
  Furthermore,
  \begin{equation*}
    \excond{S_{a}^{4}}{Z} = \sum_{I,J,K,L\ni a}\excond{W_{I}W_{J}W_{K}W_{L}}{Z}.
  \end{equation*}
  A term vanishes unless the quadruple is admissible, and it is then bounded by
  $B(I,J,K,L)$. Moreover, the quadruple is connected since $a$ belongs to the
  four sets. A quadruple is counted for at most $p$ values of $a$.
\end{proof}

\begin{theorem}[Quantitative de Jong theorem for kernels of product form]
  \label{thm-discrete-mall:dejong-product}
  Let $F$ be as above, with $\E[F^{2}] = 1$. Then
  \begin{equation*}
    d_{W}(F, \Normal(0,1)) \le \sqrt{\frac{2}{\pi}}\sum_{p=1}^{m}\var\bigl[\excond{F_{p}^{2}}{Z}\bigr]^{1/2} + C_{m}\,\mathfrak{M}^{+}(F)^{1/2},
  \end{equation*}
  where $C_{m} = 2\pi^{-1/2}m^{2}(1 + \log m) + 2\sqrt{2}\,m^{5/2}$.
\end{theorem}

\begin{proof}
  We plug Proposition~\ref{prop-discrete-mall:var-gamma-product} and
  Lemma~\ref{lem-discrete-mall:fourth-product} into
  Corollary~\ref{cor-discrete-mall:bound-multi-chaos}. Then we use
  $\sqrt{x+y}\le\sqrt{x} + \sqrt{y}$ and
  \begin{equation*}
    \sum_{p,q=1}^{m}\frac{1}{q}\le m(1+\log m), \qquad \sum_{p=1}^{m}\frac{\E[F_{p}^{2}]}{p}\le 1,
  \end{equation*}
  together with
  \begin{equation*}
    \Bigl\{\sum_{p=1}^{m}\bigl(16p\,\mathfrak{M}^{+}(F)\bigr)^{1/4}\Bigr\}^{2}\le 4m^{5/2}\,\mathfrak{M}^{+}(F)^{1/2}.
  \end{equation*}
\end{proof}

The kernels of Example~\ref{ex-discrete-mall:latent-term} are of product form,
with $\mathfrak{M}^{+}(F) = \E[c(Z)^{4}]/n$. Thus the latent term cannot be
removed from Theorem~\ref{thm-discrete-mall:dejong-product} either.

\section{Application to motif estimation}
\label{sec:appli}

We apply the previous bounds to the asymptotic normality of subhypergraph
counts in exchangeable random hypergraphs.

\subsection{Basic hypergraph definitions}
\label{sec:hypergraphs}

Hypergraphs generalise graphs and are used to model higher-order interactions
in network analysis.

\begin{definition}
  A hypergraph $G = (V, E)$ on a finite set $V = V(G)$ is given by a family $E$
  of subsets of $V$, called hyperedges. Vertices are adjacent if some hyperedge
  contains both of them; vertices contained in no hyperedge are isolated. A
  hypergraph is connected if it has no isolated vertex and the intersection
  graph of $E$ is connected.
\end{definition}
We denote by $[e]$ the set of vertices of the hyperedge $e$.
\begin{definition}
  Let $k \ge 2$. A hypergraph $G = (V, E)$ is $k$-uniform if each of its
  hyperedges has cardinality $k$, i.e.\ $E \subset \binom{V}{k}$, the
  collection of subsets of $V$ of cardinality $k$.
\end{definition}
A $2$-uniform hypergraph is a graph.
\begin{definition}
  Let $k \ge 2$ and let $G = (V, E)$ be a $k$-uniform hypergraph. A
  $k$-subhypergraph of $G$ is a $k$-uniform hypergraph $H = (V', E')$ with
  $V' \subset V$ and $E' \subset E \cap \binom{V'}{k}$.
\end{definition}
We denote by $v_H$ and $e_H$ the numbers of vertices and of hyperedges of a
hypergraph $H$. For a $3$-uniform hypergraph $G$, its shadow
$G^{(2)}$ is the set of pairs of vertices contained in some hyperedge of $G$,
and $e^{(2)}_G := |G^{(2)}|$. The automorphism group $\text{Aut}(G)$ is the
group of permutations of $V(G)$ which map $E(G)$ onto itself. A copy of
a $3$-uniform hypergraph $G$ without isolated vertices in $[n]$ is a set
$H \subset \binom{[n]}{3}$ of triples such that the hypergraph
$(\bigcup_{\alpha\in H}\alpha, H)$ is isomorphic to $G$; we then write
$H \simeq G$. There are $(n)_{v_G}/|\text{Aut}(G)|$ copies of $G$, where
$(n)_v = n(n-1)\cdots(n-v+1)$.

\subsection{Exchangeable random hypergraphs}
\label{sec:exchangeable-hypergraphs}
Random hypergraphs are natural extensions of random graphs. A large part of
the literature deals with the Erd\"os-R\'enyi model and its generalisations,
which are examples of exchangeable random hypergraphs.
\begin{definition}
  Let $k\ge2$. A $k$-uniform exchangeable random hypergraph on the vertex set
  $[n]$ is a family of $\{0,1\}$-valued random variables
  $(X_\alpha)_{\alpha \in \binom{[n]}{k}}$ such that, for every permutation
  $\sigma$ of $[n]$,
  \begin{equation*}
    (X_{\sigma(\alpha)})_{\alpha \in \binom{[n]}{k}} \overset{d}{=} (X_{\alpha})_{\alpha \in \binom{[n]}{k}},
    \qquad \sigma(\alpha) := \{\sigma(u),\, u \in \alpha\}.
  \end{equation*}
  The associated random hypergraph is $([n], E)$ with $\alpha \in E$ if and
  only if $X_\alpha = 1$.
\end{definition}
A finite exchangeable array need not be the restriction of an infinite
one, as is already the case for exchangeable sequences (draws without
replacement from an urn). If $(X_\alpha)_{\alpha\in\binom{[n]}{k}}$ is the
restriction to $[n]$ of a jointly exchangeable $\{0,1\}$-valued array indexed
by the $k$-subsets of $\N$, the Aldous--Hoover theorem (see
\cite{aldous1981representations,hoover1979relations} and
\cite[Definition~2.8]{austin_exchangeable_2008}) shows that it has the same
law as
\begin{equation}
  \label{eq-discrete-mall:aldous-hoover}
  X_{\alpha} = f\Bigl( U_{\emptyset},\ (U_{s})_{s \in \alpha},\ (U_{\beta})_{\beta \subsetneq \alpha,\ |\beta| \ge 2},\ U_{\alpha} \Bigr),
  \qquad \alpha \in \binom{[n]}{k},
\end{equation}
where $(U_{\beta})_{\beta \subseteq [n]}$ is a family of independent random
variables, uniformly distributed on $[0,1]$, and $f$ is measurable and
invariant under the permutations of $\alpha$, acting simultaneously on the
vertex sources and on the lower-order sources. We call $U_{\emptyset}$ the global
source of randomness, $(U_{s})_{s}$ the vertex sources,
$(U_{\beta})_{2 \le |\beta| < k}$ the lower-order sources and $U_{\alpha}$ the
private source of the hyperedge $\alpha$. Set
\begin{equation}
  \label{eq-discrete-mall:latent-sigma-field}
  Z := \Bigl( U_{\emptyset},\ (U_{s})_{s \in [n]},\ (U_{\beta})_{\beta \subsetneq [n],\ 2 \le |\beta| < k} \Bigr).
\end{equation}
Given $Z$, the only randomness left in $X_{\alpha}$ is its private source
$U_{\alpha}$, and the $(U_{\alpha})_{\alpha}$ are independent: the array
$(X_{\alpha})_{\alpha \in \binom{[n]}{k}}$ is therefore a family of
conditionally independent random variables given $Z$, as in
Section~\ref{sec:structure}, with $A = \binom{[n]}{k}$. We
call $Z$ the latent variable. In the sequel, we only consider arrays given by
\eqref{eq-discrete-mall:aldous-hoover}; for each $n$, the models below are of
this form.

\begin{example}[Erd\"os--R\'enyi hypergraph]
  $\mathbb{G}^{(k)}(n, p_n)$ corresponds to an $f$ depending on the private
  source only, $X_{\alpha} = \ind{U_{\alpha} \le p_n}$. The latent variable is
  trivial and the hyperedges are independent.
\end{example}

\begin{example}[Hypergraph stochastic block model]
  Let $c : [0,1] \to \llbracket 1, r \rrbracket$ assign a community to each
  vertex and let $(p_{i_{1} \ldots i_{k}})$ be a symmetric array of
  parameters in $[0,1]$. Take $X_{\alpha} = \ind{U_{\alpha} \le p_{c(U_{s}),\, s \in \alpha}}$.
  The latent variable is the community assignment, $Z = (U_{s})_{s \in [n]}$;
  given $Z$ the hyperedges are independent, but not identically distributed.
\end{example}

\begin{example}[Layered model $\mathbb{T}^{(3)}(n, q_n, p_n)$]
  \label{ex-discrete-mall:layered}
  Take $k = 3$ and
  \begin{equation}
    \label{eq-discrete-mall:layered-model}
    X_{\alpha} = \ind{U_{\alpha} \le p_n} \prod_{\beta \subset \alpha,\ |\beta| = 2} \hat{X}_{\beta},
    \qquad \hat{X}_{\beta} := \ind{U_{\beta} \le q_n}.
  \end{equation}
  A triple carries a hyperedge only if the three pairs it spans are all edges
  of the latent Erd\"os--R\'enyi graph $\mathbb{G}(n, q_n)$, and then
  independently with probability $p_n$. Since $X$ depends on the
  lower-order sources only through the latent graph, we take as latent
  variable $Z = (\hat{X}_{\beta})_{\beta \in \binom{[n]}{2}}$ instead of
  \eqref{eq-discrete-mall:latent-sigma-field}, with the
  same conditional laws for the $X_\alpha$. Then
  \begin{equation*}
    \excond{X_{\alpha}}{Z} = p_n \ind{\alpha^{(2)} \subset \mathbb{G}(n, q_n)} ,
  \end{equation*}
  where $\alpha^{(2)}$ is the set of the three pairs contained in $\alpha$.
\end{example}

For $q_n = 1$ the latent graph is complete, so that
$\mathbb{T}^{(3)}(n, 1, p_n)$ and $\mathbb{G}^{(3)}(n, p_n)$ have the same law;
for $q_n < 1$ the two models differ in many ways, as pointed out in
\cite[Example 23.11]{lovasz2012large}. In the rest
of this section, the index set $A$ of Section~\ref{sec:structure}
is the set $\binom{[n]}{3}$ of triples.

\subsection{Motif estimation: the conditionally centred count}
\label{sec:motif-estimation}

Let $G$ be a $3$-uniform hypergraph without isolated vertices. The numbers of
copies of $G$ in $\mathbb{G}^{(3)}(n, p_n)$ and in
$\mathbb{T}^{(3)}(n, q_n, p_n)$ are
\begin{equation}
  M_G =  \sum_{\substack{H \subset \binom{[n]}{3} \\H \simeq G}} \prod_{\alpha \in H} \xi_\alpha
  \quad\text{and}\quad
  N_G = \sum_{\substack{H \subset \binom{[n]}{3} \\H \simeq G}} \prod_{\alpha \in H} X_\alpha,
\end{equation}
where $\xi_\alpha = \ind{U_\alpha\le p_n}$ and $X_\alpha$ is given by
\eqref{eq-discrete-mall:layered-model}. We denote by
\begin{equation*}
  \bar M_G = \frac{M_G - \E[M_G]}{\sqrt{\var[M_G]}},\qquad
  \bar N_G = \frac{N_G - \E[N_G]}{\sqrt{\var[N_G]}},
\end{equation*}
and
\begin{equation*}
  \tilde N_G = \frac{N_G - \excond{N_G}{Z}}{\sigma_n},\qquad
  \sigma_n^2 := \var\bigl[N_G - \excond{N_G}{Z}\bigr],
\end{equation*}
the counts centred at their means and normalised by their standard
deviations, and the count of
$\mathbb{T}^{(3)}(n, q_n, p_n)$ centred at its conditional mean. The
statistic $\bar N_G$ is treated in Section~\ref{sec:fully-centred}; here we
deal with $\tilde N_G$, and with $\bar M_G$, which is the special case of a
trivial latent variable.

Write $X_\alpha = \excond{X_\alpha}{Z} + Y_\alpha$, with
$\excond{X_\alpha}{Z} = p_n\ind{\alpha^{(2)}\subset Z}$; here and below,
$\alpha^{(2)}\subset Z$ means that the three pairs of $\alpha$ are edges of the
latent graph. Note that $Y_\alpha = 0$ when $\alpha^{(2)}\not\subset Z$.
Expanding the products,
\begin{equation}
  \label{eq-discrete-mall:N_G-decomposition}
  N_G - \excond{N_G}{Z} = \sum_{j=1}^{e_G} \pi_j(N_G), \qquad
  \pi_j(N_G) = \sum_{\substack{J \subset \binom{[n]}{3}\\ |J| = j}} W_J,\qquad
  W_J = c_J \prod_{\alpha \in J} Y_\alpha,
\end{equation}
where
\begin{equation}
  \label{eq-discrete-mall:motif-count-condindep-hoeffding-term}
  c_J := \sum_{\substack{H \simeq G \\ H \supseteq J}} p_n^{|H|-|J|}\, \ind{(H \setminus J)^{(2)} \subset Z}
\end{equation}
is $\sigma(Z)$-measurable. Since the $Y_\alpha$ are conditionally independent
and conditionally centred, $W_J \in \mathfrak{C}_j$, and
\eqref{eq-discrete-mall:N_G-decomposition} is the chaos decomposition of
$N_G - \excond{N_G}{Z}$. In particular
$\sigma_n^2 = \sum_J \E[W_J^2]$. For $\mathbb{G}^{(3)}(n,p_n)$, the same
decomposition holds with $Z$ trivial, $Y_\alpha = \xi_\alpha - p_n$ and
$c_J = p_n^{e_G - |J|}\,|\{H \simeq G : H \supseteq J\}|$.

The kernels $W_J$ are of product form \eqref{eq-discrete-mall:product-form}:
the $Y_\alpha = (\xi_\alpha - p_n)\ind{\alpha^{(2)}\subset Z}$ are
conditionally independent, conditionally centred and bounded by $1$, with
$\excond{Y_\alpha^2}{Z} = p_n(1-p_n)\ind{\alpha^{(2)}\subset Z}$, and the
$c_J$ are bounded. The rate is expressed in terms of the following
quantities. For $U\subseteq V(G)$, let $e(U)$ be the number of hyperedges of
$G$ contained in $U$ and $e^{(2)}(U)$ the number of pairs of the shadow
$G^{(2)}$ contained in $U$, and set
\begin{equation}
  \label{eq-discrete-mall:phi-xi}
  \Phi_n(G) := \min_{U:\ e(U)\ge1} n^{|U|}p_n^{e(U)}q_n^{e^{(2)}(U)},\qquad
  \Xi_n(G) := \min_{U:\ e^{(2)}(U)\ge1} n^{|U|}q_n^{e^{(2)}(U)},
\end{equation}
and $\mu_n := n^{v_G}p_n^{e_G}q_n^{e^{(2)}_G}$. Below, $K$ denotes a constant
which only depends on $G$ and may change from line to line, and the sums run
over copies of $G$ in $[n]$.

\begin{lemma}
  \label{lem-motif:extension}
  Let $R$ be a set of at most $r$ vertices of $[n]$, $R^{(3)}$ a set of
  triples of elements of $R$ and $R^{(2)}$ a set of pairs of elements of
  $R$. For each copy $H$
  of $G$, let $k_H$ and $s_H$ be integers with
  $0\le k_H\le |H\cap R^{(3)}|$ and $0\le s_H\le|H^{(2)}\cap R^{(2)}|$. Then,
  with $K_r$ depending only on $G$ and $r$,
  \begin{align*}
    \sum_{H:\ H\cap R^{(3)}\neq\emptyset} p_n^{e_G-k_H}q_n^{e^{(2)}_G-s_H}   & \le K_r\,\frac{\mu_n}{\Phi_n(G)}, \\
    \sum_{H:\ H^{(2)}\cap R^{(2)}\neq\emptyset} p_n^{e_G}q_n^{e^{(2)}_G-s_H} & \le K_r\,\frac{\mu_n}{\Xi_n(G)}.
  \end{align*}
\end{lemma}

\begin{proof}
  Group the copies $H$ according to the set $U\subseteq V(G)$ of vertices of
  $G$ sent into $R$ by an isomorphism, and to the images of these vertices.
  There are at most $2^{v_G}r^{v_G}$ classes, each containing at most
  $n^{v_G-|U|}$ copies. In a class, the triples of $H\cap R^{(3)}$ and the
  pairs of $H^{(2)}\cap R^{(2)}$ are contained in the image of $U$, so that
  $k_H\le e(U)$ and $s_H\le e^{(2)}(U)$. Since $p_n, q_n\le1$, it follows that
  \begin{equation*}
    p_n^{e_G-k_H}\le p_n^{e_G-e(U)} \quad\text{and}\quad q_n^{e^{(2)}_G-s_H}\le q_n^{e^{(2)}_G-e^{(2)}(U)}.
  \end{equation*}
  In the first sum, $e(U)\ge1$ and a class contributes at most
  \begin{equation*}
    \frac{\mu_n}{n^{|U|}p_n^{e(U)}q_n^{e^{(2)}(U)}}\le\frac{\mu_n}{\Phi_n(G)}.
  \end{equation*}
  In the second one, $e^{(2)}(U)\ge1$ and a class contributes at most
  \begin{equation*}
    \frac{\mu_n}{n^{|U|}q_n^{e^{(2)}(U)}}\le\frac{\mu_n}{\Xi_n(G)}.
  \end{equation*}
\end{proof}

\begin{lemma}
  \label{lem-motif:variance}
  For $n\ge 2v_G$, $p_n\in(0,1)$ and $q_n\in(0,1]$,
  \begin{equation*}
    K^{-1}(1-p_n)\,\frac{\mu_n^2}{\Phi_n(G)} \le \sigma_n^2 \le K(1-p_n)\,\frac{\mu_n^2}{\Phi_n(G)}.
  \end{equation*}
  In particular $\sigma_n^2 > 0$, so that $\tilde N_G$ is well defined.
\end{lemma}

\begin{proof}
  By conditional orthogonality,
  \begin{equation*}
    \sigma_n^2 = \sum_J \bigl(p_n(1-p_n)\bigr)^{|J|}\,\E\bigl[c_J^2\ind{J^{(2)}\subset Z}\bigr],
  \end{equation*}
  and
  \begin{equation*}
    c_J^2\ind{J^{(2)}\subset Z} = \sum_{H,H'\supseteq J}p_n^{2e_G-2|J|}\ind{H^{(2)}\cup H'^{(2)}\subset Z}.
  \end{equation*}
  Let $k = |H\cap H'|$ be the number of common hyperedges of $H$ and $H'$.
  Then
  \begin{equation*}
    \sum_{\emptyset\neq J\subseteq H\cap H'}\Bigl(\frac{1-p_n}{p_n}\Bigr)^{|J|} = p_n^{-k}-1.
  \end{equation*}
  Thus, exchanging the sums,
  \begin{equation}
    \label{eq-discrete-mall:sigma-pairs}
    \sigma_n^2 = \sum_{(H,H'):\ k\ge1} p_n^{2e_G}\bigl(p_n^{-k}-1\bigr)q_n^{|H^{(2)}\cup H'^{(2)}|},
  \end{equation}
  where
  \begin{equation*}
    (1-p_n)p_n^{-k}\le p_n^{-k}-1\le k(1-p_n)p_n^{-k}.
  \end{equation*}
  For the upper bound, the sum over $H$ of $p_n^{e_G}q_n^{e^{(2)}_G}$ is at
  most $\mu_n$. The sum over $H'$ is bounded by
  Lemma~\ref{lem-motif:extension} with
  \begin{equation*}
    R = V(H),\quad R^{(3)} = H,\quad R^{(2)} = H^{(2)},\quad k_{H'} = k \quad\text{and}\quad s_{H'} = |H^{(2)}\cap H'^{(2)}|.
  \end{equation*}
  For the lower bound, let $U^*$ achieve $\Phi_n(G)$. For each copy $H$,
  consider the copies $H'$ which coincide with $H$ on the image of $U^*$ and
  whose other vertices lie outside $V(H)$. There are at least
  $K^{-1}n^{v_G-|U^*|}$ of them, and they satisfy
  \begin{equation*}
    k = e(U^*)\ge1 \quad\text{and}\quad |H^{(2)}\cap H'^{(2)}| = e^{(2)}(U^*).
  \end{equation*}
  Each such pair contributes at least
  \begin{equation*}
    (1-p_n)\,p_n^{2e_G-e(U^*)}\,q_n^{2e^{(2)}_G-e^{(2)}(U^*)}
  \end{equation*}
  to \eqref{eq-discrete-mall:sigma-pairs}. Since there are at least
  $K^{-1}n^{v_G}$ copies $H$, the lower bound follows.
\end{proof}

\begin{theorem}
  \label{thm-motif_estimation:dejong-quant-hypergraph}
  Let $G$ be a $3$-uniform hypergraph without isolated vertices. For
  $n\ge2v_G$, $p_n\in(0,1)$ and $q_n\in(0,1]$,
  \begin{equation*}
    d_W\bigl(\tilde N_G,\Normal(0,1)\bigr) \le K\Biggl(\frac{1}{\sqrt{(1-p_n)\Phi_n(G)}} + \sqrt{\frac{1-q_n}{\Xi_n(G)}}\Biggr).
  \end{equation*}
\end{theorem}

\begin{proof}
  We apply Theorem~\ref{thm-discrete-mall:dejong-product} to
  $\tilde N_G$, whose kernels are the $W_J/\sigma_n$, with $m = e_G$.

  We first bound the connected sum. Since
  \begin{equation*}
    c_J\ind{J^{(2)}\subset Z} = \sum_{H\supseteq J}p_n^{e_G-|J|}\ind{H^{(2)}\subset Z},
  \end{equation*}
  one has
  \begin{equation*}
    \mathfrak{M}^+(\tilde N_G) \le \sigma_n^{-4}\sum p_n^{4e_G - \sum_i|J_i| + |\cup_iJ_i|}(1-p_n)^{|\cup_iJ_i|}\,q_n^{|\cup_iH_i^{(2)}|},
  \end{equation*}
  the sum running over $(J_1,\dots,J_4)\in\mathcal{Q}$ and copies
  $H_i\supseteq J_i$. Since $J_i\subseteq H_i$,
  \begin{equation*}
    \sum_i|J_i|-\Bigl|\bigcup_iJ_i\Bigr| \le \sum_{\alpha\in\cup_iJ_i}\bigl(\#\{i:\alpha\in H_i\}-1\bigr) \le 4e_G - \Bigl|\bigcup_iH_i\Bigr|,
  \end{equation*}
  so that the power of $p_n$ is at least $|\cup_iH_i|$. Moreover,
  $(1-p_n)^{|\cup_iJ_i|}\le 1-p_n$. Since $(J_i)$ is connected, so is the
  graph on $\{H_1,\dots,H_4\}$ joining two copies with a common hyperedge.
  Order the copies so that each one shares a hyperedge with an earlier one,
  and let $k_i$ and $s_i$ be the numbers of hyperedges and of shadow pairs of
  $H_i$ which belong to $H_1\cup\dots\cup H_{i-1}$. Then
  \begin{equation*}
    p_n^{|\cup_i H_i|}q_n^{|\cup_i H_i^{(2)}|} = \prod_{i=1}^4 p_n^{e_G-k_i}q_n^{e^{(2)}_G-s_i},
  \end{equation*}
  with $k_i\ge1$ for $i\ge2$. Summing over $H_1$, then over $H_2$, $H_3$,
  $H_4$ with Lemma~\ref{lem-motif:extension}, and accounting for the
  $2^{4e_G}$ choices of $J_i\subseteq H_i$ and the $4!$ orders, we get
  \begin{equation*}
    \mathfrak{M}^+(\tilde N_G) \le K(1-p_n)\,\sigma_n^{-4}\,\mu_n\Bigl(\frac{\mu_n}{\Phi_n(G)}\Bigr)^3 \le \frac{K}{(1-p_n)\Phi_n(G)}
  \end{equation*}
  by Lemma~\ref{lem-motif:variance}.

  We now bound the latent term. As in the proof of
  Lemma~\ref{lem-motif:variance},
  \begin{equation*}
    \excond{\pi_j(N_G)^2}{Z} = \sum_{(H,H')} a_j(H,H')\,\ind{H^{(2)}\cup H'^{(2)}\subset Z},
  \end{equation*}
  where, with $k = |H\cap H'|$, $a_j(H,H') = 0$ if $k<j$ and otherwise
  \begin{equation*}
    a_j(H,H') = \binom{k}{j}p_n^{2e_G-j}(1-p_n)^{j}\le 2^{e_G}(1-p_n)p_n^{2e_G-k}.
  \end{equation*}
  For sets $A$, $B$ of pairs,
  \begin{equation*}
    \cov(\ind{A\subset Z},\ind{B\subset Z}) = q_n^{|A\cup B|}\bigl(1-q_n^{|A\cap B|}\bigr).
  \end{equation*}
  This covariance vanishes if $A\cap B = \emptyset$. Otherwise, it is at most
  \begin{equation*}
    2e^{(2)}_G\,(1-q_n)\,q_n^{|A\cup B|}.
  \end{equation*}
  Hence $\var[\excond{\pi_j(N_G)^2}{Z}]$ is at most $K(1-p_n)^2(1-q_n)$ times
  the sum of
  \begin{equation*}
    p_n^{4e_G-k_1-k_2}q_n^{|H_1^{(2)}\cup H_1'^{(2)}\cup H_2^{(2)}\cup H_2'^{(2)}|}
  \end{equation*}
  over the copies $(H_1,H_1',H_2,H_2')$ such that
  \begin{equation*}
    k_1 = |H_1\cap H_1'|\ge1,\quad k_2 = |H_2\cap H_2'|\ge1 \quad\text{and}\quad \bigl(H_1^{(2)}\cup H_1'^{(2)}\bigr)\cap\bigl(H_2^{(2)}\cup H_2'^{(2)}\bigr)\neq\emptyset.
  \end{equation*}
  Up to a factor $2$, we may assume that $H_2^{(2)}$ meets
  $H_1^{(2)}\cup H_1'^{(2)}$. We sum over $H_1$, which gives at most $\mu_n$.
  Then we sum over $H_1'$ with the first bound of
  Lemma~\ref{lem-motif:extension} and $R = V(H_1)$. Next, we sum over $H_2$
  with the second bound and $R = V(H_1)\cup V(H_1')$. Finally, we sum over
  $H_2'$ with the first bound and $R = V(H_1)\cup V(H_1')\cup V(H_2)$. Here
  $k_2$ is at most the number of hyperedges of $H_2'$ in
  $R^{(3)} = H_1\cup H_1'\cup H_2$. With Lemma~\ref{lem-motif:variance}, this
  gives
  \begin{equation*}
    \frac{\var[\excond{\pi_j(N_G)^2}{Z}]}{\sigma_n^4} \le K\,\frac{(1-p_n)^2(1-q_n)\,\mu_n^4\,\Phi_n(G)^{-2}\,\Xi_n(G)^{-1}}{(1-p_n)^2\mu_n^4\,\Phi_n(G)^{-2}} = K\,\frac{1-q_n}{\Xi_n(G)}.
  \end{equation*}
  Theorem~\ref{thm-discrete-mall:dejong-product} concludes the proof.
\end{proof}

When $q_n = 1$, the latent graph is complete, $\tilde N_G = \bar M_G$, and
$\Phi_n(G)$ is the minimum of $n^{v_H}p_n^{e_H}$ over the subhypergraphs $H$
of $G$ with at least one hyperedge, since for a fixed vertex set the induced
subhypergraph has the most hyperedges. We recover, in the Wasserstein
distance, the rate of \cite[Theorem~3.3]{michalczuk2024normal}, including the
factor $1-p_n$; the corresponding central limit theorem goes back to
\cite{de1996central}.

\begin{corollary}
  \label{cor-motif:erdos-renyi}
  Let $G$ be a $3$-uniform hypergraph without isolated vertices. For
  $n\ge2v_G$ and $p_n\in(0,1)$,
  \begin{equation}
    \label{eq-discrete-mall:hypergraph-condindep-threshold}
    d_{W}(\bar M_G, \Normal(0, 1)) \le K \Bigl( (1-p_n)\min_{\substack{H \subseteq G \\ e_H \ge 1}} n^{v_H} p_n^{e_H} \Bigr)^{-1/2}.
  \end{equation}
\end{corollary}

For $p_n\le c<1$, the same rate holds for the Kolmogorov distance by
\cite[Theorem~3.3]{michalczuk2024normal}.

\begin{example}[One hyperedge]
  \label{ex-motif:one-hyperedge}
  If $G$ is a single hyperedge, then $\Phi_n(G) = n^3p_nq_n^3$ and
  $\Xi_n(G) = \min(n^2q_n, n^3q_n^3)$, so that
  \begin{equation*}
    d_W(\tilde N_G,\Normal(0,1))\le K\Biggl(\frac{1}{\sqrt{n^3p_n(1-p_n)q_n^3}} + \sqrt{1-q_n}\,\Bigl(\frac{1}{\sqrt{n^2q_n}} + \frac{1}{(nq_n)^{3/2}}\Bigr)\Biggr).
  \end{equation*}
  The latent term of this bound is of the order of $\sqrt{\var[T]}/\E[T]$,
  where $T$ is the number of triangles of the latent graph, and the condition
  $nq_n\to\infty$ is necessary for the asymptotic normality of $\tilde N_G$:
  $\tilde N_G = 0$ on $\{T = 0\}$, and by the Harris inequality
  $\mathbb{P}(T=0)\ge(1-q_n^3)^{\binom n3}$, which is bounded away from $0$
  when $nq_n$ is bounded.
\end{example}

\subsection{The fully centred count}
\label{sec:fully-centred}

We now consider the count $N_G$ of $\mathbb{T}^{(3)}(n,q_n,p_n)$ centred at its
mean. By the tower property, its variance splits orthogonally into a
hyperedge part and a latent part,
\begin{equation}
  \label{eq-discrete-mall:variance-split}
  \var[N_G] = \E\bigl[\var(N_G \mid Z)\bigr] + \var\bigl(\excond{N_G}{Z}\bigr),
\end{equation}
where $\excond{N_G}{Z} = p_n^{e_G}\,\bigl|\{H \simeq G : H^{(2)} \subset \mathbb{G}(n,q_n)\}\bigr|$
is a subgraph count of the latent graph. The normal approximation of
$\bar N_G$ follows from the results of \cite{michalczuk2024normal} for
non-uniform random hypergraphs with independent edges, through the following
identification. Let $\tilde G$ be the
non-uniform hypergraph with vertex set $V(G)$ and edge set $E(G) \cup G^{(2)}$:
it carries the $3$-edges of $G$ and the $2$-edges of its shadow. Let
$\mathbb{H}(n; q_n, p_n)$ be the random hypergraph on $[n]$ in which every pair
is an edge with probability $q_n$ and every triple is an edge with probability
$p_n$, all independently, and let $Z^n_{\tilde G}$ be the number of its
subhypergraphs isomorphic to $\tilde G$.

\begin{lemma}
  \label{lem-motif:identification}
  With $\xi_\alpha = \ind{U_\alpha \le p_n}$ and $\hat X_\beta$ as in
  \eqref{eq-discrete-mall:layered-model}, the random hypergraph with $3$-edges
  $\{\alpha : \xi_\alpha = 1\}$ and $2$-edges $\{\beta : \hat X_\beta = 1\}$ has
  the law of $\mathbb{H}(n; q_n, p_n)$, and on this probability space
  $N_G = Z^n_{\tilde G}$ almost surely. In particular $\bar N_G$ has the law of
  the standardised count of $\tilde G$ in $\mathbb{H}(n;q_n,p_n)$.
\end{lemma}

\begin{proof}
  The variables $(\xi_\alpha)_\alpha$ and $(\hat X_\beta)_\beta$ are
  independent Bernoulli variables with parameters $p_n$ and $q_n$, which gives
  the first claim. For a copy $H$ of $G$, since the $\hat X_\beta$ are
  $\{0,1\}$-valued, repeated factors collapse and
  \begin{equation*}
    \prod_{\alpha \in H} X_\alpha = \prod_{\alpha \in H}\xi_\alpha\prod_{\alpha\in H}\prod_{\beta\subset\alpha}\hat X_\beta
    = \prod_{\alpha \in H}\xi_\alpha\prod_{\beta \in H^{(2)}}\hat X_\beta ,
  \end{equation*}
  which is the indicator that the copy $H \cup H^{(2)}$ of $\tilde G$ is present.
  The map $H \mapsto H \cup H^{(2)}$ is a bijection from the copies of $G$ onto
  the copies of $\tilde G$: an isomorphism maps the shadow of $G$ onto the
  shadow of its image, and the inverse map keeps the $3$-edges.
\end{proof}

For a subhypergraph $F$ of $\tilde G$, denote by $v_F$, $e_3(F)$ and $e_2(F)$
its numbers of vertices, $3$-edges and $2$-edges, and set
\begin{equation}
  \label{eq-discrete-mall:psi}
  \Psi_n(G) := \min_{\substack{F \subseteq \tilde G \\ e_3(F) + e_2(F) \ge 1}} n^{v_F}\, p_n^{e_3(F)}\, q_n^{e_2(F)} .
\end{equation}
For a fixed vertex set, the minimum is attained by the induced subhypergraph,
so that $\Psi_n(G)$ is also the minimum over the sets $W\subseteq V(G)$
spanning at least one edge of $\tilde G$ of
$n^{|W|}p_n^{|E(G)[W]|}q_n^{|G^{(2)}[W]|}$.

\begin{theorem}
  \label{thm-discrete-mall:hypergraph-condindep-threshold-hoeffding}
  Let $G$ be a $3$-uniform hypergraph without isolated vertices and assume that
  $0 < p_n, q_n \le c < 1$. Then
  \begin{equation}
    d_{W}(\bar N_G, \Normal(0, 1)) + d_{K}(\bar N_G, \Normal(0, 1)) \lesssim \Psi_n(G)^{-1/2},
  \end{equation}
  where $d_K$ is the Kolmogorov distance and the implicit constant only depends
  on $G$ and $c$. Moreover, $\bar N_G$ converges in distribution to
  $\Normal(0,1)$ if and only if $\Psi_n(G) \to \infty$.
\end{theorem}

\begin{proof}
  By Lemma~\ref{lem-motif:identification}, $\bar N_G$ has the law of the
  standardised count of $\tilde G$ in the non-uniform random hypergraph
  $\mathbb{H}(n;q_n,p_n)$, in which each edge is present independently with a
  probability depending only on its size. The bound is then
  \cite[Theorem~3.3, (3.4)]{michalczuk2024normal}, which applies since all
  edge probabilities are at most $c$. Indeed, for $F\subseteq\tilde G$, the
  quantity of \cite{michalczuk2024normal}
  \begin{equation*}
    P_F = \prod_{e\in E(F)}p_{|e|} \quad\text{is equal to}\quad p_n^{e_3(F)}q_n^{e_2(F)},
  \end{equation*}
  so that the minimum in (3.4) is $\Psi_n(G)$. Furthermore, the probabilities
  of the edges of sizes other than $2$ and $3$ do not enter. The
  characterisation is \cite[Theorem~3.2]{michalczuk2024normal}. Since
  $p_n, q_n \le c$, the second condition there, which involves
  $n^{|e|}(1-p_{|e|})$, is automatically satisfied.
\end{proof}

The sets $W$ with $|E(G)[W]|\ge1$ give $\Phi_n(G)$, and the sets with
$|E(G)[W]| = 0$ appear in $\Xi_n(G)$; the sets with $|E(G)[W]|\ge1$ appear in
$\Xi_n(G)$ with the larger value $n^{|W|}q_n^{|G^{(2)}[W]|}$. Hence
$\Psi_n(G) = \min(\Phi_n(G),\Xi_n(G))$, and
Theorem~\ref{thm-motif_estimation:dejong-quant-hypergraph} gives
\begin{equation*}
  d_W\bigl(\tilde N_G,\Normal(0,1)\bigr)\le K\bigl((1-p_n)\Psi_n(G)\bigr)^{-1/2}.
\end{equation*}
Under the assumptions of
Theorem~\ref{thm-discrete-mall:hypergraph-condindep-threshold-hoeffding}, the
condition $\Psi_n(G)\to\infty$, which characterises the asymptotic normality
of $\bar N_G$, thus implies that of $\tilde N_G$, at least at the same rate.

\section*{Funding}
The second author was supported in part by LabEx Math\'ematiques Hadamard.

\providecommand{\bysame}{\leavevmode\hbox to3em{\hrulefill}\thinspace}
\providecommand{\MR}{\relax\ifhmode\unskip\space\fi MR }
\providecommand{\MRhref}[2]{%
  \href{http://www.ams.org/mathscinet-getitem?mr=#1}{#2}
}
\providecommand{\href}[2]{#2}


\begin{thebibliography}{10}

\bibitem{aldous1981representations}
David~J. Aldous, \emph{Representations for partially exchangeable arrays of
  random variables}, Journal of Multivariate Analysis \textbf{11} (1981),
  no.~4, 581--598.

\bibitem{austin_exchangeable_2008}
Tim Austin, \emph{On exchangeable random variables and the statistics of large
  graphs and hypergraphs}, Probability Surveys \textbf{5} (2008), no.~none,
  80--145.

\bibitem{barbour1989central}
Andrew~D Barbour, Micha{\l} Karo{\'n}ski, and Andrzej Ruci{\'n}ski, \emph{A
  central limit theorem for decomposable random variables with applications to
  random graphs}, Journal of Combinatorial Theory, Series B \textbf{47} (1989),
  no.~2, 125--145.

\bibitem{chen2010normal}
Louis~HY Chen, Larry Goldstein, and Qi-Man Shao, \emph{Normal approximation by
  stein’s method}, Springer Science \& Business Media, 2010.

\bibitem{de1996central}
Peter De~Jong, \emph{{A Central Limit Theorem with Applications to Random
  Hypergraphs}}, Random Structures and Algorithms \textbf{8} (1996), no.~2,
  105--120.

\bibitem{Decreusefond2022}
Laurent Decreusefond, \emph{Selected topics in malliavin calculus: Chaos,
  divergence and so much more}, Springer Nature, 2022.

\bibitem{decreusefond:hal-01565240}
Laurent Decreusefond and H{\'e}l{\`e}ne Halconruy, \emph{{Malliavin and
  Dirichlet structures for independent random variables}}, {Stochastic
  Processes and their Applications} \textbf{129} (2019), no.~8, 2611--2653.

\bibitem{dobler2023normal}
Christian D{\"o}bler, \emph{Normal approximation via non-linear exchangeable
  pairs}, Latin American Journal of Probability and Mathematical Statistics
  \textbf{20} (2023), 167.

\bibitem{dobler2026malliavin}
\bysame, \emph{On the {M}alliavin calculus on product spaces and an infinite de
  {J}ong theorem}, arXiv:2609.15725, 2026.

\bibitem{dobler2017quantitative}
Christian D{\"o}bler and Giovanni Peccati, \emph{{Quantitative de Jong theorems
  in any dimension}}, Electronic Journal of Probability \textbf{22} (2017).

\bibitem{duerinckx2021size}
Mitia Duerinckx, \emph{{On the size of chaos via {G}lauber calculus in the
  classical mean-field dynamics}}, Communications in Mathematical Physics
  (2021), 1--41.

\bibitem{dung2018poisson}
Nguyen~Tien Dung, \emph{{Poisson and normal approximations for the measurable
  functions of independent random variables}}, arXiv preprint arXiv:1807.10925
  (2018).

\bibitem{hoover1979relations}
Douglas~N. Hoover, \emph{Relations on probability spaces and arrays of random
  variables}, Tech. report, Institute for Advanced Study, Princeton, 1979.

\bibitem{kallenberg1997foundations}
Olav Kallenberg, \emph{Foundations of modern probability}, Springer, 1997.

\bibitem{krokowski2017discrete}
Kai Krokowski, Anselm Reichenbachs, and Christoph Th{\"a}le, \emph{Discrete
  {Malliavin--Stein} method: {Berry}--{Esseen} bounds for random graphs and
  percolation}, The Annals of Probability \textbf{45} (2017), no.~2,
  1071--1109.

\bibitem{lovasz2012large}
L{\'a}szl{\'o} Lov{\'a}sz, \emph{Large networks and graph limits}, vol.~60,
  American Mathematical Soc., 2012.

\bibitem{mcdiarmid1989method}
Colin McDiarmid, \emph{On the method of bounded differences}, Surveys in
  combinatorics \textbf{141} (1989), no.~1, 148--188.

\bibitem{michalczuk2024normal}
Wojciech Michalczuk, Miko{\l}aj Nieradko, and Grzegorz Serafin, \emph{Normal
  approximation for subgraph count in random hypergraphs}, arXiv preprint
  arXiv:2408.06112 (2024).

\bibitem{mossel2010noise}
Elchanan Mossel, Ryan O'Donnell, and Krzysztof Oleszkiewicz, \emph{Noise
  stability of functions with low influences: invariance and optimality},
  Annals of Mathematics \textbf{171} (2010), no.~1, 295--341.

\bibitem{NourdinNormalApproximationsMalliavin2012}
Ivan Nourdin and Giovanni Peccati, \emph{Normal approximations with malliavin
  calculus: from stein's method to universality}, Cambridge Tracts in
  Mathematics, no. 192, Cambridge University Press, 2012.

\bibitem{privault2008stochastic}
Nicolas Privault, \emph{{Stochastic analysis of Bernoulli processes}},
  Probability Surveys \textbf{5} (2008), 435--483.

\bibitem{privault2020normal}
Nicolas Privault and Grzegorz Serafin, \emph{Normal approximation for sums of
  discrete {$ U $}-statistics -- {A}pplication to {Kolmogorov} bounds in random
  subgraph counting}, Bernoulli \textbf{26} (2020), no.~1, 587--615.

\bibitem{rucinski1988small}
Andrzej Ruci{\'n}ski, \emph{When are small subgraphs of a random graph normally
  distributed?}, Probability Theory and Related Fields \textbf{78} (1988),
  no.~1, 1--10.

\end{thebibliography}

\end{document}